\documentclass[draft]{amsart}

\usepackage{amssymb}
\usepackage{url}
\usepackage{amscd}

\theoremstyle{plain}
\newtheorem{theorem}{Theorem}

\newtheorem{corollary}[theorem]{Corollary}
\newtheorem{lemma}[theorem]{Lemma}
\newtheorem{proposition}[theorem]{Proposition}

\theoremstyle{remark}
\newtheorem{remark}[theorem]{Remark}
\newtheorem{remark and definition}[theorem]{Remark and Definition}

\newtheorem{example}[theorem]{Example}

\newtheorem{fact}[theorem]{Fact}

\newcommand{\R}{\mathbb{R}}
\newcommand{\Q}{\mathbb{Q}}

\newcommand{\Sf}{\mathbb{S}}

\newcommand{\Les}{\mathbb{L}}

\newcommand{\Oes}{\mathbb{O}}

\newcommand{\spa}{\mbox{span}}

\newcommand{\rank}{\mbox{rank }}

\newcommand{\grad}{\mbox{grad}}

\newcommand{\ii}{isometric immersion }
\newcommand{\iis}{isometric immersions }
\newcommand{\su}{submanifold }
\newcommand{\sus}{submanifolds }

\newcommand{\sff}{second fundamental form }

\newcommand\fall{\;\;\;\mbox{for all}\;\;}
\newcommand\an{\;\;\mbox{and}\;\;}

\def\V{{\mathcal V}}

\def\Fes{{\mathcal F}}

\def\H{{\mathcal H}}

\def\D{{\mathcal D}}
\def\<{\langle}

\def\>{\rangle}
\def\a{\alpha}
\def\va{\varphi}

\def\o{\omega}

\def\d{\partial}
\def\bea{\begin{eqnarray*} }
\def\eea{\end{eqnarray*} }
\def\be{\begin{equation} }
\def\ee{\end{equation} }

\def\proof{\noindent{\it Proof: }}
\def\qed{\ifhmode\unskip\nobreak\fi\ifmmode\ifinner
\else\hskip5 pt \fi\fi\hbox{\hskip5 pt \vrule width4 pt
height6 pt  depth1.5 pt \hskip 1pt }}

\begin{document}


\title[Isometric immersions of warped products]{Isometric immersions in codimension two of
\\ warped products into space forms
}
\author{Marcos Dajczer and Ruy Tojeiro}

{\renewcommand{\baselinestretch}{1}


\address{IMPA -- Estrada Dona
Castorina 110\\
22460-320-Rio de Janeiro\\Brazil
}\email{marcos@impa.br}\vspace{2ex}
\address{
Universidade Federal de S\~{a}o Carlos\\
13565-905-S\~{a}o Carlos\\Brazil} \email{tojeiro@dm.ufscar.br}


\subjclass{53B25}

\begin{abstract}
We provide a local classification of \iis $f\colon\,L^p\times_\rho
M^n\to\Q_c^{\,p+n+k}$ in codimensions $k=1, 2$ of warped products
of Riemannian manifolds into space forms, under the assumptions
that $n\geq k+1$ and that $N^{p+n}=L^p\times_\rho M^n$ has no
points with the same constant sectional curvature $c$ as the
ambient space form.

\end{abstract}

\maketitle


\newcommand\sfrac[2]{{#1/#2}}

\newcommand\cont{\operatorname{cont}}
\newcommand\diff{\operatorname{diff}}


\section{Introduction}
A basic decomposition theorem due to N\"olker \cite{no} states
that  an isometric immersion $f\colon\, N^{p+n}=L^p\times_\rho
M^n\to\Q_c^{\,\ell}$ of a warped product of connected Riemannian
manifolds with warping function $\rho\in C^{\infty}(L^p)$ into a
complete simply-connected space form of constant sectional
curvature $c$ is a {\it warped product of isometric immersions\/}
(see \cite{no} or Section $1$ for the precise definition of this
concept) whenever its \sff $\a\colon\,TN\times TN\to T^\perp N$
satisfies
$$
\a(X,V)=0\fall X\in TL\an V\in TM.
$$
This generalizes a well-known result for isometric immersions of
Riemannian products into Euclidean space due to Moore \cite{mo} as
well as its extension by Molzan \cite{mol} for nonflat ambient
space forms; see also \cite{sc}.

   It is a natural problem to understand the possible
cases in which the \ii $f$ may fail, locally or globally, to be a
warped product of isometric immersions. In high codimensions the
warped product structure of the manifold does not seem to place
enough restrictions on the \ii in order to make possible a
complete classification in either case.  Even in the much more
restrictive situation of Riemannian products, a successful local
analysis has only been carried out in the case in which the
codimension is two and the first factor is one-dimensional; see
\cite{bdt}.  This was used therein to characterize \iis $f\colon\,
L^p\times M^n\to\R^{p+n+2}$ of complete Riemannian products none
of whose factors is everywhere flat (see Remark~\ref{bdt} below),
carrying through the global results previously obtained  by Moore
\cite{mo} and Alexander-Maltz \cite{am}. Based on earlier work due
to  Moore, such \sus were shown  in \cite{am} to split as products
of hypersurfaces under the global assumption that they do not
carry an Euclidean strip.

  The main goal of this  paper is to provide a local classification
of \iis $f\colon\,L^p\times_\rho M^n\to\Q_c^{\,p+n+k}$ in
codimensions $k=1, 2$,  under the assumptions that $n\geq k+1$ and
that $N^{p+n}=L^p\times_\rho M^n$ has no points with the same
constant sectional curvature $c$ as the ambient space form.  In
the case of codimension $k=1$, we prove that $f$ must be a warped
product of isometric immersions. In codimension $k=2$ we show that
only two other possibilities may arise. Namely, either $f$ is a
composition of a warped product of isometric immersions into
$\Q_c^{\,p+n+1}$ with a local \ii of $\Q_c^{\,p+n+1}$ into
$\Q_c^{\,p+n+2}$ or $N^{p+n}$ is a Riemannian manifold of a
special type that admits a second decomposition as a warped
product with respect to which $f$ splits as a warped product of
isometric immersions. We leave the precise statement for Section
$3$, where we also state its corresponding version for the case of
Riemannian products.

   We give examples showing that the restriction on the
dimension of $M^n$ is necessary. As for the hypothesis that
$N^{p+n}$ has no points with constant sectional curvature~$c$, we
observe that for manifolds with constant sectional curvature the
assumption that they be warped products places no further
restrictions on them since any such manifold can be realized as a
warped product in many possible ways; see the discussion on warped
product representations of space forms in Section $1$. Therefore,
for Riemannian manifolds with constant sectional curvature our
problem reduces to classifying {\it all\/} \iis in codimension two
of such manifolds into space forms. In this regard, recall that if
$f\colon\,U\subset \Q_{\tilde{c}}^{\,n}\to \Q_c^{\,n+2}$, $n\geq
4$, is an isometric immersion, then $\tilde{c}\geq c$ and, for
$\tilde{c}> c$, away from the set of umbilical points the
immersion must be locally a composition of the umbilical inclusion
of $\Q_{\tilde{c}}^{\,n}$ into $\Q_c^{\,n+1}$ with a local \ii of
$\Q_c^{\,n+1}$ into $\Q_c^{\,n+2}$ (see \cite{er}, \cite{he} and
\cite{dt1}).
 In fact, the latter statement can be derived from our main theorem, but in that result we exclude from our analysis
the case of local \iis of $\Q_c^{\,n}$ into $\Q_c^{\,n+2}$ with
the same constant sectional curvature. A local description of
these \iis when $c=0$ was given in \cite{df}.

As striking applications of our main result, we obtain that if
$L^p$ is a Riemannian manifold no open subset of which can be
isometrically immersed into $\Q_c^{\,p+1}$, then any \ii
$f\colon\,L^p\times_{\rho}M^n\to \Q_c^{\,p+n+2}$, $n\geq 3$, is
either a cylindrical submanifold of Euclidean space or a
rotational submanifold. Moreover, if $L^p$ cannot be locally
isometrically immersed in $\Q_c^{\,p+2}$, then
$L^p\times_{\rho}M^n$ cannot be locally isometrically immersed in
$\Q_c^{\,p+n+2}$ either for whatever Riemannian manifold $M^n$ of
dimension $n\geq 3$ and warping function~$\rho$.\vspace{1ex}

\noindent {\it Acknowledgment.\/}  We are greatly indebted to the
referee whose many suggestions and comments have decisively
contributed for significant improvements in both the presentation
and the mathematical content of this article.


\section{Preliminaries}

In this section we establish our notation and state some basic
facts on warped products of Riemannian manifolds and their
isometric immersions into the standard real space forms.
\vspace{1ex}

 Given a a vector bundle $E$ over a Riemannian manifold $N$,
we denote by $\Gamma(E)$ the set of all locally defined smooth
sections of $E$.  If $N=L\times M$ is a product manifold,  we
denote by $\H$ and $\V$ the {\it horizontal\/} and {\it
vertical\/} subbundles of $TN$, that is, the distributions on $N$
correspondent to the product foliations determined by $L$ and $M$,
respectively. Elements of $\Gamma(\H)$ will always be denoted by
the letters $X,Y,Z$, whereas those in $\Gamma(\V)$ by the letters
$U,V,W$. The same applies to individual tangent vectors. A vector
field $X\in \Gamma(\H)$ (resp., $V\in \Gamma(\V)$) is said to be
the {\em lift\/} of a vector field $\tilde{X}\in \Gamma(TL)$
(resp., $\tilde{V}\in \Gamma(TM)$) if ${\pi_L}_*X=\tilde{X}\circ
\pi_L$ (resp., ${\pi_M}_*V=\tilde{V}\circ \pi_M$), where
\mbox{$\pi_L\colon\,L\times M\to L$} (resp., $\pi_M\colon\,L\times
M\to M$) is the canonical projection onto $L$ (resp., $M$). We
denote the set of all lifts of vector  fields in  $L$ (resp., $M$)
by ${\mathcal L}(L)$  \mbox{(resp., ${\mathcal L}(M)$)}, and  always
denote vector fields in $L$ and $M$ with a tilde and use the same
letters without the tilde to represent their  lifts to $N$.

    If $L$ and $M$ are Riemannian manifolds with Riemannian metrics
$\<\;,\;\>_L$ and $\<\;,\;\>_M$, respectively, the {\em warped
product\/} $N=L\times _{\rho}M$ with {\em warping function\/}
$\rho\in C^{\infty}(L)$ is the product manifold $L\times M$
endowed with the {\em warped product metric\/}
$$
\<\;\,,\,\;\>={\pi_L}^*\<\;\,,\,\;\>_L+(\rho\circ
\pi_L)^2{\pi_M}^*\<\,\;,\,\;\>_M.
$$
We always assume $N$ to be connected. The Levi-Civita connections
of $N$, $L$ and $M$ are related by (cf.\ \cite{on}) \be \nabla_X
Y\;\;\mbox{is the lift of} \;\;\nabla^L_{\tilde{X}}
\tilde{Y},\label{eq:c1}\ee \be\nabla_X V =\nabla_V
X=-\<X,\eta\>V,\label{eq:c2}\ee \be(\nabla_V W)_{\V}\;\; \mbox{is
the lift of} \;\; \nabla^M_{\tilde{V}} \tilde{W},\label{eq:c3} \ee
\be (\nabla_V W)_{\H}=\<V,W\>\eta,\label{eq:c4} \ee where $X,Y\in
{\mathcal L}(L)$, $V,W\in {\mathcal L}(M)$ and
$\eta=-\,\grad\,(\log\rho\circ \pi_{L})$. Here and in the sequel,
writing a vector field with a vector subbundle as a subscript
indicates  taking the section of that vector subbundle obtained by
orthogonally projecting the vector field pointwise onto the
corresponding fiber of the subbundle.  Observe that the formula
$\nabla_X V =-\<X,\eta\>V$ (resp., $\nabla_V X=-\<X,\eta\>V$) is
tensorial in $X$ (resp., $V$), hence it also holds for horizontal
(resp., vertical)  vector fields that are not necessarily lifts.
On the other hand, it characterizes vertical (resp., horizontal)
vector fields that are lifts.

   Recall that a vector subbundle $E$ of $TN$ is called
{\em totally geodesic\/}  or {\em autoparallel\/} if $\nabla_X
Y\in \Gamma(E)$ for all $X,Y\in \Gamma(E)$. It is called {\em
totally umbilical\/} if there exists a vector field $\eta\in
\Gamma(E^\perp)$ such that $(\nabla_X Y)_{E^\perp}=\<X,Y\>\eta$
for all $X,Y\in \Gamma(E)$. If, in addition, the so called {\em
mean curvature normal\/} $\eta$ of $E$ satisfies $(\nabla_X
\eta)_{E^\perp}=0$ for all $X\in \Gamma(E)$, then $E$ is said to
be {\em  spherical\/}.  A totally umbilical vector subbundle $E$
of $TN$ is automatically integrable, and its leaves are totally
umbilical submanifolds of $N$. If $E$ is totally geodesic or
spherical then the leaves are totally geodesic or spherical
submanifolds of $N$, respectively. By a {\it spherical\/}
submanifold we mean a totally umbilical submanifold whose mean
curvature vector field is parallel in the normal connection.

    It follows from (\ref{eq:c1}) and (\ref{eq:c4}), respectively,
    that ${\mathcal H}$ is totally geodesic and that $\V$ is totally umbilical
    with mean curvature normal $\eta=-\,\grad\,(\log\rho\circ \pi_{L})$.
Moreover,  since $\eta$ is a gradient vector field and $\H$ is
totally geodesic we have
$$
\<\nabla_V \eta, X\>=\<\nabla_X \eta, V\>=0
$$
for all $X\in \Gamma(\H)$ and $V\in \Gamma(\V)$, and hence $\V$ is
spherical. The following extension due to Hiepko  of the
well-known decomposition theorem of de Rham shows that these
properties characterize warped products.

\begin{theorem}{\em (\cite{hi})}\label{thm:hi}
Let $N$ be a Riemannian manifold and let $TN=\H\oplus \V$ be an
orthogonal decomposition into nontrivial vector subbundles such
that $\H$ is  totally geodesic and $\V$ is spherical. Then, for
every point $z_0\in N$ there exist an isometry $\Psi$  of a warped
product $L\times _{\rho}M$ onto a neighborhood of $z_0$ in $N$
such that $\Psi(L\times \{x\})$ and $\Psi(\{y\}\times M)$  are
integral manifolds of $\H$ and $\V$, respectively, for all $y\in
L$ and $x\in M$. Moreover, if $N$ is simply connected and complete
then the isometry $\Psi$ can be taken  onto all of $N$.
\end{theorem}

  Given a warped product  $N=L\times _{\rho}M$,
the lift of the curvature tensor ${}^M\!R$ of $M$ to $N$ is the
tensor whose value at $E_1, E_2, E_3\in T_z N$ is the unique
vector in ${\mathcal V}_z$ that projects to
${}^M\!R({\pi_M}_*E_1,{\pi_M}_*E_2){\pi_M}_*E_3$ in
$T_{\pi_M(z)}M$.  The lift of the curvature tensor ${}^LR$ of $L$
is similarly defined. Then the curvature tensors of $L$, $M$ and
$N$ are related by
\begin{eqnarray}
\label{eq:i}R(X,Y)Z\!\!\!&=&\!\!\!{}^LR(X,Y)Z,\\
\label{eq:iii}R(X,Y)V\!\!\!&=&\!\!\!R(V,W)X=0,\\
\label{eq:iv}R(X,U)V\!\!\!&=&
\!\!\!\<U,V\>(\nabla_X\eta-\<\eta,X\>\eta),\\
\label{eq:v}R(V,W)U\!\!\!&=&\!\!\!{}^M\!R(V,W)U
-\|\eta\|^2(\<W,U\>V-\<V,U\>W).
\end{eqnarray}
 Since
$\nabla_X\eta-\<\eta,X\>\eta\in {\mathcal H}$ because ${\mathcal
H}$ is totally geodesic,  all the information of (\ref{eq:iv}) is
contained in \be\label{eq:ivb} \<R(X,V)W,Y\>
=\<V,W\>\<\nabla_X\eta-\<X,\eta\>\eta, Y\>. \ee

  The starting point for the proof of the main results of this paper is the observation
  that the curvature relations (\ref{eq:i}) to  (\ref{eq:v}) impose several restrictions
  on the second fundamental form
$\alpha\colon\;TN\times TN\to T^\perp N$ of an  isometric
immersion $f\colon\,N\to\Q_c^{\,\ell}$ when combined with the
Gauss equation for $f$.

\begin{proposition}\label{prop:symm} Let $f\colon\,L\times _{\rho}M\to\Q_c^{\,\ell}$ be an isometric immersion of a warped product. Then the curvature-like tensor
\begin{eqnarray}
C(E_1,E_2,E_3,E_4):\!\!\!&=&\!\!\!\<R(E_1,E_2)E_3,E_4\> -
c\<(E_1\wedge E_2)E_3,E_4\> \nonumber \vspace{1ex}\\
\!\!\!&=&\!\!\!\<\a(E_1,E_4),\a(E_2,E_3)\>
-\<\a(E_1,E_3),\a(E_2,E_4)\>
\end{eqnarray}
satisfies
\begin{eqnarray}
C(X,V,W,Y) \!\!&=&\!\!\<V,W\>\<\nabla_X\eta-\<X,\eta\>\eta-cX,Y\>,
\label{eq:ge11}\\
C(X,Y,V,Z)\,\!\!&=&\!\!0,\label{eq:ge22}\\
C(X,Y,V,W)\!\!&=&\!\!0,\label{eq:ge33}\\
C(X,U,V,W)\!\!&=&\!\!0\label{eq:ge44}.
\end{eqnarray}
\end{proposition}

 We now introduce  the notion of a warped product of
isometric immersions  into $\Q_c^{\,\ell}$ which plays a
fundamental  role in this paper. This relies on the  warped
product representations of $\Q_c^{\,\ell}$, that is, isometries of
warped products onto open subsets of  $\Q_c^{\,\ell}$. All such
isometries were described by N\"olker  for warped products with
arbitrarily many factors; see \cite{no} for details.  In
particular, any isometry  of a  warped product with two factors
onto an open subset of $\Q_c^{\,\ell}$ arises as a restriction of
an explicitly constructible isometry
$$
\Psi\colon\,V^{\ell-m}(\subset \Q_c^{\,\ell-m}) \times
_{\sigma}\Q_{\tilde{c}}^{\,m}\to\Q_c^{\,\ell}
$$
onto an open dense subset of $\Q_c^{\,\ell}$, where
$\Q_{\tilde{c}}^{\,m}$ is a complete spherical submanifold of
$\Q_c^{\,\ell}$ and $V^{\ell-m}$ is an open subset of the unique
totally  geodesic submanifold $\Q_c^{\,\ell-m}$ of $\Q_c^{\,\ell}$
(of constant sectional curvature $c$ if $\ell-m\ge 2$) whose
tangent space at some point $\bar{z}\in \Q_{\tilde{c}}^{\,m}$ is
the orthogonal complement of the tangent space of
$\Q_{\tilde{c}}^{\,m}$ at $\bar{z}$. The  isometry $\Psi$ is, in
fact, completely determined by the choice  of
$\Q_{\tilde{c}}^{\,m}$ and of a point $\bar{z}\in
\Q_{\tilde{c}}^{\,m}$, and it is called the
 {\em warped product representation\/} of $\Q_c^{\,\ell}$
determined by  $(\bar{z},\Q_{\tilde{c}}^{\,m})$. If $c\neq 0$, we
consider the standard model of $\Q_c^{\,\ell}$ as a complete
spherical submanifold of $\Oes^{\,\ell+1}$, where
$\Oes^{\,\ell+1}$ denotes either Euclidean space $\R^{\ell+1}$ if
$c>0$ or Lorentzian space $\Les^{\ell+1}$ if $c<0$. Then, for
$c\neq 0$ the warping function $\sigma$ is the restriction to
$V^{\ell-m}$ of the height function $z\mapsto \<z,a\>$ in
$\Oes^{\,\ell+1}$, where $-a$ is the mean curvature vector of
$\Q_{\tilde{c}}^{\,m}$ in $\Oes^{\,\ell+1}$ at $\bar{z}$.
Similarly, if $c=0$ then $\sigma(z)= 1+\<z-\bar{z},a\>$, where
$-a$ is the mean curvature vector of $\Q_{\tilde{c}}^{\,m}$ in
$\Q_c^{\,\ell}=\R^\ell$ at $\bar{z}$. In every case
$\<a,a\>=\tilde{c}$.\vspace{1ex}

\noindent {\bf Definition.} Let
$\Psi\colon\,{V^{\ell-m}}\times_\sigma \Q_{\tilde c}^{\,m} \to
\Q_c^{\,\ell}$ be a warped product representation, let
\mbox{$h_1\colon\,L\to V^{\ell-m}$} and $h_2\colon\,M\to
\Q_{\tilde c}^{\,m}$ be isometric immersions, and let
$\rho=\sigma\circ h_1$. Then the \ii $f=\Psi\circ(h_1\times
h_2)\colon\,N= L\times_{\rho} M\to\Q_c^{\,\ell}$ is called the
{\it warped product of  the \iis $h_1$ and $h_2$ determined by
$\Psi$\/}.

\bigskip
\vspace{1ex}

\begin{picture}(150,84)\hspace*{-10ex}
\put(100,30){$L\times_{\rho=\sigma\circ h_1} M$}
\put(110,55){$h_1$} \put(165,32){\vector(1,0){135}}
\put(305,30){$\Q_c^{\,\ell}$} \put(158,55){$h_2$}
\put(205,45){$\circlearrowright$}
\put(180,20){$f=\Psi\circ(h_1\times h_2)$}
\put(105,42){\vector(0,1){30}} \put(152,42){\vector(0,1){30}}
\put(102,80){$V^{\ell-m}\times_\sigma \Q_{\tilde c}^{\,m}$}
\put(175,80){\vector(3,-1){125}} \put(235,65){$\Psi$}
\end{picture}
\begin{example}\label{ex:1}  If $N=L\times_\rho M$
is not a Riemannian product and $h_2$ is an isometry, then $f$ is
called a {\it rotational submanifold\/} with profile $h_1$. This
means that $V^{\ell-m}$ is a half-space of a totally geodesic
submanifold $\Q_c^{\,\ell-m}\subset \Q_c^{\,\ell}$ bounded by a
totally geodesic submanifold $\Q_c^{\,\ell-m-1}$ and $f(N)$ is the
submanifold of $\Q_c^{\,\ell}$ generated by the action on $h_1(L)$
of the subgroup of isometries of $\Q_c^{\,\ell}$ that leave
$\Q_c^{\,\ell-m-1}$ invariant.
\end{example}

\begin{example}\label{ex:2}  If $N=L\times_\rho M$
is not a Riemannian product and $h_1\colon\,L\to V^{\ell-m}$ is a
local isometry then, for $c=0$, we have that $f(N)$ is contained
in the product of an  Euclidean factor $\R^{\ell-m-1}$ with a cone
in $\R^{m+1}$ over $h_2$. If $c\neq 0$, then $f(N)$ is the union
of open subsets of the totally geodesic \sus of $\Q_c^{\,\ell}$
through the points of $h_2(M)\subset \Q_{\tilde c}^{\,m}$ whose
tangent spaces at the points of $h_2(M)$ are the normal spaces of
$\Q_{\tilde c}^{\,m}$ in $\Q_c^{\,\ell}$.
\end{example}

  Notice that any warped product of isometric immersions
in codimension one must be as in one of the preceding examples. In
codimension two only a third possibility arises, namely, the case
in which both $h_1$ and $h_2$ are hypersurfaces.

    Important special cases of warped products of
isometric immersions arise as follows. Let  $\Q_{c_1}^{\,\ell_1}$
and $\Q_{c_2}^{\,\ell_2}$ be complete spherical submanifolds of
$\Q_c^{\,\ell}$ through a fixed point $\bar{z}\in \Q_c^{\,\ell}$
whose tangent spaces  at $\bar{z}$ are orthogonal and whose mean
curvature vectors $\psi_1$ and $\psi_2$ at $\bar{z}$ satisfy
$\<\psi_1,\psi_2\>=-c$ and $\psi_1$ (resp., $\psi_2$) is
orthogonal to $T_{\bar{z}}\Q_{c_2}^{\,\ell_2}$ (resp.,
$T_{\bar{z}}\Q_{c_1}^{\,\ell_1}$).  Let
$\Psi\colon\,V^{\ell-\ell_2}\times _{\sigma}\Q_{c_2}^{\,\ell_2}
\to\Q_c^{\,\ell}$ be the warped product representation of
$\Q_c^{\,\ell}$ determined by $(\bar{z},\Q_{c_2}^{\,\ell_2})$.
Then $\Q_{c_1}^{\,\ell_1}\subset V^{\ell-\ell_2}$ and $\sigma\circ
i=1$, where $i\colon\, \Q_{c_1}^{\,\ell_1}\to V^{\ell-\ell_2}$ is
the inclusion map. The warped product $\Psi\circ (i\times
id)\colon\,\Q_{c_1}^{\,\ell_1}\times \Q_{c_2}^{\,\ell_2}\to
\Q_c^{\,\ell}$ of the inclusion and the identity map is an
isometric embedding called the isometric embedding of the
Riemannian product $\Q_{c_1}^{\,\ell_1}\times
\Q_{c_2}^{\,\ell_2}$ into $\Q_c^{\,\ell}$ as an {\em extrinsic
Riemannian product\/}.\vspace{1ex}

The structure of the \sff of a warped product of \iis is described
in the following result.

\begin{proposition} \label{prop:sffs}
Let $N=L\times_\rho M$ and $\bar{N}=\bar{L}\times_{\bar{\rho}}
\bar{M}$ be warped product manifolds, and let
$F\colon\,L\to\bar{L}$ and $G\colon\, M\to\bar{M}$ be \iis with
$\rho=\bar{\rho}\circ F$. Then $f=F\times G\colon\,N\to\bar{N}$ is
an \ii and at $z=(y,x)\in N$ we have
\begin{itemize}
\item[$(i)$] ${\pi_{\bar{L}}}_*f_*T_zN =F_* T_yL$,
$\;{\pi_{\bar{L}}}_*T^\perp_z N=T_y^\perp L$,
$\;{\pi_{\bar{M}}}_*f_*T_zN =G_* T_xM$,
$\;{\pi_{\bar{M}}}_*T^\perp_z N=T_x^\perp M$.\vspace{1ex}
\item[$(ii)$]  $(\grad\,
\bar{\rho}\,(F(y)))_{T_y^\perp L} =\grad\, \bar{\rho}\,(F(y))-F_*
\grad\, \rho\,(y)$.\vspace{1ex}
\item[$(iii)$] The second fundamental forms of $F$,
$G$ and $f$ are related by \be\label{eq:sff2}
{\pi_{\bar{M}}}_*\a^f(E_1,E_2) =\a^G({\pi_M}_*E_1, {\pi_M}_*E_2),
\ee \be\label{eq:sff1}\begin{array}{l}
{\pi_{\bar{L}}}_*\a^f(E_1,E_2)=\a^F({\pi_L}_*E_1,
{\pi_L}_*E_2)-\vspace{1ex}\\\hspace{17.5ex}\rho\,(y)\<{\pi_M}_*E_1,
{\pi_M}_*E_2\> (\grad\, \bar{\rho}\,(F(y)))_{T_y^\perp L}.
\end{array}\ee
\end{itemize}
\end{proposition}

   Given an \ii $f\colon\,N\to\bar{N}$, a normal vector
$\zeta\in T_z^\perp N$ is called a {\it principal curvature normal
vector\/} at $z$ if the subspace
$$
\Delta_\zeta(z)=\{T\in T_zN:\a(T,E) =\<T,E\>\zeta \mbox{ for
all}\;E\in T_zN\}
$$
is nontrivial. In this case $\Delta_\zeta(z)$ is called the
eigenspace correspondent to $\zeta$. If $\zeta=0$ then
$\Delta(z):=\Delta_0(z)$ is called the {\it relative nullity
subspace\/} of $f$ at $z$.

\begin{corollary} \label{le:triv} Let $f$
be an \ii as in Proposition \ref{prop:sffs}. Then at $z=(y,x)\in
N$ we have
\begin{itemize}
\item[$(i)$]  $\V_z\subset \Delta_\zeta(z)$ for a principal curvature normal vector
$\zeta\in T_z^\perp N$ if and only if  $G$ is umbilical at $x$
with mean curvature vector ${\pi_{\bar{M}}}_*\zeta$ and
${\pi_{\bar{L}}}_*\zeta =-\rho^{-1}(\grad\,
\bar{\rho}\,)_{T_{y}^\perp L}$. In particular, we have that
$\V_z\subset \Delta(z)$ if and only if  $G$ is totally geodesic at
$x$ and~$(\grad\, \bar{\rho}\,)_{T_{y}^\perp L}=~0$.
\item[$(ii)$] $\H_z\subset
\Delta(z)$  if and only if $F$ is totally geodesic at $y$.
\end{itemize}
\end{corollary}

       Given a vector $a\neq 0$ in either
$\R^{\ell}$ or $\Oes^{\,\ell+1}$, according as $c=0$ or $c\neq 0$,
let $U$ be the vector field on $\Q_c^{\,\ell}$ defined by
$U_z=a-c\<a,z\>z$ and let $\Fes^a$ be the $1$-dimensional totally
geodesic distribution generated by $U$ on the open dense subset
$W^a=\{z\in \Q_c^{\,\ell}: U_z\neq 0\}$. Notice that
$\Q_c^{\,\ell}\setminus W^a$ is empty for $c=0$ as well as for
$(c<0, \<a,a\>\geq 0)$, it contains one point for $(c<0, \<a,a\><
0)$ and two points for $c>0$. Observe also that for $c=0$ (resp.,
$c\neq 0$) the vector field  $U$ is the gradient of the function
$\sigma\colon\;\Q_c^{\,\ell}\to \R$ given by
$\sigma(z)=1+\<z-\bar{z},a\>$ for a fixed $\bar{z} \in
\Q_c^{\,\ell}$ (resp., $\sigma(z)=\<z,a\>$), which was used in the
definition of a warped product representation of $\Q_c^{\,\ell}$.
We say that an isometric immersion $g\colon\, L^p\to
\Q_c^{\,\ell}$ is  {\em cylindrical with respect to $a$ \/} if
$g(L)\subset W^a$ and $\Fes^a$ is everywhere tangent to $g(L)$, or
equivalently, if $U_{g(y)}=\grad\; \sigma(g(y))$ is nonzero and
tangent to $g(L)$ for any $y\in L$. The last assertion in
Corollary~\ref{le:triv}-$(i)$ yields the following.

\begin{corollary} \label{cor:cyl} Let $f=\Psi\circ (F\times G)\colon\,N=L\times_\rho M\to \Q_c^{\,\ell}$ be a warped product of isometric immersions, where $\Psi\colon\,{V^{\ell-m}}\times_\sigma \Q_{\tilde c}^{\,m}
\to \Q_c^{\,\ell}$ is a warped product representation  determined
by $(\bar{z},\Q_{\tilde{c}}^{m})$. If $G$ is totally geodesic and
$F$ is cylindrical with respect to the mean curvature vector $-a$
of $\Q_{\tilde{c}}^{m}$ at $\bar{z}$ in either $\R^{\ell}$ or
$\Oes^{\,\ell+1}$, according as $c=0$ or $c\neq 0$, then the
vertical subbundle of $TN$ is contained in the relative nullity
subbundle of $f$. Conversely, if the the vertical subbundle of
$TN$ is contained in the relative nullity subbundle of $f$ then
$G$ is totally geodesic and $F|_U$ is cylindrical with respect to
the mean curvature vector $-a$ of $\Q_{\tilde{c}}^{m}$ at
$\bar{z}$, where $U$ is the open subset of $L$ where $\grad\;\rho$
does not vanish.\end{corollary}

  If $f\colon\,L\times_{\rho} M\to \Q_c^{\,\ell}$ is a
warped product of isometric immersions, then it follows from
Proposition \ref{prop:sffs}-$(ii)$ that at any point $z\in L\times
M$ its \sff satisfies \be\label{eq:condA} \a(X,V)=0\fall X\in
\H_z\an V\in \V_z. \ee The following theorem due to N\"olker
states that the converse is also true. Recall that the {\em
spherical hull\/} of an \ii $G\colon\,M\to \Q_c^{\,\ell}$ is the
complete spherical  \su of least dimension that contains $G(M)$.

\begin{theorem}{\em (\cite{no})}\label{thm:no}
Let $f\colon\,L\times_\rho M\to\Q_c^{\,\ell}$ be an \ii of a
warped product  whose \sff satisfies condition (\ref{eq:condA})
everywhere. For a fixed point $(\bar{y},\bar{x})\in L\times_\rho
M$ with $\rho(\bar{y})=1$,  let $F\colon\,L\to \Q_c^{\,\ell}$ and
$G\colon\,M\to\Q_c^{\,\ell}$ be given by $F(y)=f(y,\bar{x})$ and
$G(x)=f(\bar{y},x)$, and let $\Q_{\tilde c}^{\,m}$ be the
spherical hull of $G$.   Then $(f(\bar{y},\bar{x}), \Q_{\tilde
c}^{\,m})$ determines a warped product representation
$\Psi\colon\,V^{\ell-m} \times _{\sigma}\Q_{\tilde{c}}^{\,m} \to
\Q_c^{\,\ell}$ such that $F(L)\subset V^{\ell-m}$ and $f=\Psi
\circ (F\times G)$, where in the last equation $F$ and $G$ are
regarded as maps into $V^{\ell-m}$ and $\Q_{\tilde{c}}^{\,m}$,
respectively.
\end{theorem}

 The preceding theorem is also valid for \iis of warped products
with arbitrarily  many factors (see \cite{no}). It contains as a
particular case the following result due to Molzan (cf.\ Corollary
$17$ of \cite{no}), which is an extension  to nonflat ambient
space forms of the main lemma  in \cite{mo}.

\begin{corollary} {\em (\cite{mol})}\label{thm:molzan}
Let $f\colon\,L\times M\to\Q_c^{\,\ell}$ be an isometric immersion
of a Riemannian product  whose second fundamental form satisfies
condition (\ref{eq:condA}) everywhere. For a fixed point
$(\bar{y},\bar{x})\in L\times M$  define $F\colon\,L\to
\Q_c^{\,\ell}$ and $G\colon\,M\to \Q_c^{\,\ell}$ by
$F(y)=f(y,\bar{x})$ and $G(x)=f(\bar{y},x)$, and denote by
$\Q_{c_1}^{\,\ell_1}$ and $\Q_{c_2}^{\,\ell_2}$ the spherical
hulls of $F(L)$ and $G(M)$, respectively. Then $F$ and $G$ are
isometric immersions and there exists an isometric embedding
$\Phi\colon\,\Q_{c_1}^{\,\ell_1}\times\Q_{c_2}^{\,\ell_2}
\to\Q_c^{\,\ell}$ as an extrinsic  Riemannian product such that
\mbox{$f=\Phi\circ (F\times G)$}, where in the last equation $F$
and $G$ are regarded as maps into $\Q_{c_1}^{\,\ell_1}$ and
$\Q_{c_2}^{\,\ell_2}$, respectively.
\end{corollary}

    In applying Theorem \ref{thm:no} one must often be able to determine the dimension of the spherical hull of $G$. In the remaining of this section we develop a tool for computing that dimension.

    Given an \ii $g\colon\,M^n\to \Q_c^{\,\ell}$, a subbundle $\tilde{{\mathcal Z}}$
    of the normal bundle of $g$ is called {\em umbilical\/} if there exists
    $\theta\in \Gamma(\tilde{{\mathcal Z}})$ such that
$$
(\alpha^g(E_1,E_2))_{\tilde{{\mathcal Z}}}=\<E_1,E_2\>\theta
$$
for all $E_1, E_2\in \Gamma(TM)$. We say that $\theta$ is the {\em
principal curvature normal\/}
 of $\tilde{{\mathcal Z}}$. If $n\geq 2$ and the subbundle $\tilde{{\mathcal Z}}$ is parallel in the normal connection,
  then the Codazzi equations of $g$ imply that the vector field $\theta$ is also parallel in
  the normal connection. In particular, it has constant length. If $g(M^{n})$ is contained in a complete spherical
submanifold $\Q_{\tilde{c}}^m$ of $\Q_c^{\,\ell}$ with dimension
$m$ and constant sectional curvature $\tilde{c}$, then the
pulled-back subbundle $\tilde{{\mathcal
Z}}=g^*T^\perp\Q_{\tilde{c}}^m$, where $T^\perp\Q_{\tilde{c}}^m$
is the normal bundle of $\Q_{\tilde{c}}^m$ in $\Q_c^{\,\ell}$, is
an umbilical parallel subbundle of $T^\perp M$ of rank $\ell-m$.
Conversely, we have the following  result due to Yau.

\begin{proposition} {\em (\cite{ya})}
\label{prop:yau} Let  $g\colon\,M^n\to \Q_c^{\,\ell}$, $n\geq 2$,
be an isometric immersion. Assume that there exists an umbilical
parallel subbundle $\tilde{{\mathcal Z}}$ of $T^\perp M$ with
principal curvature normal $\theta$ and rank $\ell-m$. Then there
exists a complete spherical submanifold $\Q_{\tilde{c}}^m$ of
$\Q_c^{\,\ell}$ with dimension $m$ and constant sectional
curvature $\tilde{c}=c+\|\theta\|^2$ such that $g(M^{n})\subset
\Q_{\tilde{c}}^m$.
\end{proposition}

   As a consequence, the dimension of the spherical hull
of an isometric immersion can be characterized as follows.

\begin{corollary}\label{prop:shull} Let  $g\colon\,M^n\to \Q_c^{\,\ell}$, $n\geq 2$,
be an isometric immersion. Then the dimension of the spherical
hull of $g$ is $m$ if and only if $\ell -m$ is the maximal rank of
an umbilical parallel subbundle $\tilde{{\mathcal Z}}$ of $T^\perp M$.
Moreover, the spherical hull of $g$ has constant sectional
curvature $\tilde{c}=c+\|\theta\|^2$, where $\theta$ is the
principal curvature normal of $\tilde{{\mathcal Z}}$.
\end{corollary}

\begin{corollary}\label{cor:shull} Let $f\colon\,N^{p+n}=L^p\times_\rho M^n\to\Q_c^{\,\ell}$, $n\geq 2$, be an \ii
of a warped product  whose \sff satisfies  (\ref{eq:condA})
everywhere. Given $\bar{y}\in L$ with $\rho(\bar{y})=1$, let
$G\colon\,M\to\Q_c^{\,\ell}$ be defined by $G=f\circ i_{\bar{y}}$,
where $i_{\bar{y}}\colon\,M^n\to N^{p+n}$ given by
$i_{\bar{y}}(x)=(\bar{y},x)$  is the (isometric) inclusion of
$M^n$ into $N^{p+n}$ as a leaf of the vertical subbundle
${\mathcal V}$. Then the spherical hull of $G$ has dimension
$m=\ell-p-k$, where $k$ is the maximal rank of a parallel
subbundle ${\mathcal Z}$ of $i_{\bar{y}}^*T^\perp N$ such that
\be\label{eq:psub} \alpha^f({i_{\bar{y}}}_*V,
{i_{\bar{y}}}_*W)_{{\mathcal Z}}=\<V,W\>\theta \ee  for some
$\theta\in \Gamma({\mathcal Z})$ and for all $V,W\in \Gamma(TM)$.
If ${\mathcal Z}$ is such a subbundle, then $\theta\in
\Gamma({\mathcal Z})$ is parallel, hence has constant length.
Moreover, the spherical hull of $G$ has constant sectional
curvature  $c+\|\theta\|^2+\|\grad\; log\rho(\bar{y})\|^2$.
\end{corollary}
\proof     The normal bundle of $i_{\bar{y}}$ is
$i_{\bar{y}}^*{\mathcal H}$, where $\H$ is the horizontal subbundle of
$TN$, hence the normal bundle of $G$ splits as
$$
T_G^\perp M=i_{\bar{y}}^*T^\perp N\oplus f_*i_{\bar{y}}^*{\mathcal H}
$$
and the  second fundamental form of $G$ splits accordingly as
\be\label{for} \alpha^G(V,W)=\alpha^f({i_{\bar{y}}}_*V,
{i_{\bar{y}}}_*W)+\<V,W\>f_*(\eta\circ i_{\bar{y}}), \ee where
$\eta=-\,\grad\,(\log\rho\circ \pi_{L})$ is the mean curvature
normal of $\V$. In particular, it follows that
$f_*i_{\bar{y}}^*{\mathcal H}$ is an umbilical subbundle of $T_G^\perp
M$ with principal curvature normal $f_*(\eta\circ i_{\bar{y}})$.
Moreover, using that the \sff  of $f$ satisfies (\ref{eq:condA})
it follows that $f_*i_{\bar{y}}^*{\mathcal H}$ is parallel in the
normal connection of $G$. It is now easily seen that a subbundle
${\mathcal Z}$ of $i_{\bar{y}}^*T^\perp N$ is parallel and satisfies
(\ref{eq:psub}) if and only if ${\mathcal Z}\oplus
f_*i_{\bar{y}}^*{\mathcal H} $ is a parallel umbilical subbundle of
$T_G^\perp M$ with principal curvature normal
$\theta+f_*(\eta\circ i_{\bar{y}})$. The conclusion follows from
Corollary \ref{prop:shull}.\vspace{1.5ex}\qed

   In the sequel only the following two special cases of Corollary \ref{cor:shull}
   will be needed, in which the assumptions in part $(i)$ (resp., $(ii)$) easily imply that
   the vector subbundle ${\mathcal Z}$ equals $i_{\bar{y}}^*T^\perp N$ (resp.,
   $\{0\}$).

\begin{corollary}\label{cor:shull2} Under the
assumptions of Corollary \ref{cor:shull} we have
\begin{itemize}
\item[$(i)$] If the vertical subbundle ${\mathcal V}$ is contained in the eigendistribution correspondent to a principal curvature normal $\zeta$ of $f$, then the
spherical hull of $G$ has dimension $m=n$ and constant sectional
curvature $\tilde{c}= c+\|\zeta\circ i_{\bar{y}}\|^2+\|\grad\;
log\,\rho(\bar{y})\|^2$.
\item[$(ii)$] If there exists no local vector field $\bar{\xi}\in \Gamma(i_{\bar{y}}^*T^\perp N)$
such that $A^f_{\bar{\xi}}\circ
{i_{\bar{y}}}_*=\lambda{i_{\bar{y}}}_*$ for some $\lambda\in
C^\infty(M)$, then the spherical hull of $G$ has dimension
$m=\ell-p$ and constant sectional curvature $\tilde{c}=c+\|\grad\;
log\,\rho(\bar{y})\|^2$.
\end{itemize}
\end{corollary}

\section{The results}

 Our main result provides a complete local classification
of \iis $f\colon L^p\times_\rho M^n\to\Q_c^{\,p+n+2}$ of a  warped
product  under the assumptions that $n\geq 3$ and that
$N^{p+n}=L^p\times_\rho M^n$ is free of points with constant
sectional curvature~$c$.  Here and in the sequel it is always
assumed that $p,n\geq 1$, and only further restrictions on those
dimensions are explicitly stated.

\begin{theorem}\label{thm:main}  Assume that a warped product $N^{p+n}=L^p\times_\rho
M^n$ with $n\geq 3$ is free of points with constant sectional
curvature $c$. Then for any \ii \mbox{$f\colon
N^{p+n}\to\Q_c^{\,p+n+2}$} there exists an open dense subset of
$N^{p+n}$ each of whose points lies in an open  product
neighborhood $U=L^p_0\times M^n_0\subset L^p\times M^n$ such that
one of the following possibilities holds:
\begin{itemize}
\item[$(i)$]  $f|_U$ is a warped product of \iis
with respect to a warped product representation \mbox{$\Psi\colon
V^{p+k_1}\!\times\! _{\sigma} \Q_{\tilde{c}}^{\,n+k_2}\to
\Q_c^{\,p+n+2}$},  $k_1+k_2=2$.

\bigskip

\begin{picture}(150,84)\hspace{-15ex}
\put(100,30){$L^p_0\times_{\sigma\circ h_1} M^n_0$}
\put(112,55){$h_1$} \put(165,32){\vector(1,0){135}}
\put(305,30){$\Q_c^{\,p+n+2}$} \put(154,55){$h_2$}
\put(205,45){$\circlearrowright$}
\put(180,20){$f|_U=\Psi\circ(h_1\times h_2)$}
\put(107,42){\vector(0,1){30}} \put(149,42){\vector(0,1){30}}
\put(98,80){$V^{p+k_1}\times_\sigma \Q_{\tilde c}^{\,n+k_2}$}
\put(175,80){\vector(3,-1){125}} \put(235,65){$\Psi$}
\end{picture}
\vspace*{-2ex}
\item[$(ii)$] $f|_U$  is a composition $H\circ g$ of \iis
where $g$  is a warped~product of \iis
\mbox{$g=\Psi\circ(h_1\times h_2)$} determined by a warped product
representation $\Psi\colon\,V^{p+k_1}\times_{\sigma}
\Q_{\tilde{c}}^{\,n+k_2} \to \Q_c^{\,p+n+1}$ with $k_1+k_2=1$, and
$H\colon\, W\to \Q_c^{\,p+n+2}$ is an \ii of an open
subset~$W\supset g(U)$ of~$\,\Q_c^{\,p+n+1}$.

\bigskip
\bigskip

\begin{picture}(150,84)\hspace*{-20ex}
\put(100,30){$L^p_0\times_{\sigma\circ h_1} M^n_0$}
\put(110,55){$h_1$} \put(168,32){\vector(1,0){132}}
\put(305,30){$\Q_c^{\,p+n+2}$} \put(154,55){$h_2$}
\put(225,55){$\circlearrowright$}
\put(176,18){$f|_U=H\circ\Psi\circ(h_1\times h_2)$}
\put(106,42){\vector(0,1){30}} \put(147,42){\vector(0,1){30}}
\put(98,80){$V^{p+k_1}\times_\sigma \Q_{\tilde c}^{\,n+k_2}$}
\put(175,80){\vector(1,0){125}} \put(225,86){$\Psi$}
\put(305,77){$W\subset\Q_c^{\,p+n+1}$} \put(312,55){$H$}
\put(308,72){\vector(0,-1){30}}
\end{picture}
\vspace*{-2ex}
\item[$(iii)$] There exist open intervals $I,J\subset\R$
such that $L_0^p$,  $M_0^n$, $U$ split as
$L_0^p=L_0^{p-1}\times_{\rho_1}I$, $M_0^n
=J\times_{\rho_2}M_0^{n-1}$ and
$$
U =L_0^{p-1}\times_{\rho_1}
((I\times_{\rho_3}J)\times_{\bar{\rho}}M_0^{n-1}),
$$
where $\rho_1\in C^\infty (L_0^{p-1})$, $\rho_2\in C^\infty (J)$,
$\rho_3\in C^\infty (I)$ and $\bar{\rho}\in C^\infty (I\times J)$
satisfy
$$
\rho=(\rho_1\circ \pi_{L_0^{p-1}})(\rho_3\circ \pi_{I})
\;\;\mbox{and}\;\;\bar{\rho} =(\rho_3\circ \pi_{I})(\rho_2\circ
\pi_{J}),
$$
and there exist warped product representations
$$
\hspace*{8ex}\Psi_1\colon\,V^{p-1}\times_{\sigma_1}
\Q_{\tilde{c}}^{\,n+3}\to \Q_c^{\,p+n+2}\an
\Psi_2\colon\,W^{4}\times_{\sigma_2} \Q_{\bar{c}}^{\,n-1}\to
\Q_{\tilde{c}}^{\,n+3},
$$
an \ii $g\colon\, I\times_{\rho_3}J\to W^{4}$ and isometries
$i_1\colon\,L_0^{p-1} \to W^{\,p-1}\subset V^{\,p-1}\subset
\Q_c^{\,p-1}$ and $i_2\colon\,M^{n-1} \to W^{\,n-1}\subset
\Q_{\bar{c}}^{\,n-1}$ onto open subsets such that
$f|_U=\Psi_1\circ (i_1\times (\Psi_2\circ (g\times i_2)))$,
$\bar{\rho} =\sigma_2\circ g$ and $\rho_1=\sigma_1\circ i_1$.
Moreover, $L_0^p$ has constant sectional curvature $c$ if $p\geq
2$.

\bigskip
\bigskip

\begin{picture}(150,84)\hspace*{-15ex}
\put(115,0){$W^{\,4}\;\,\times_{\sigma_2}\, \;\Q_{\bar
c}^{\,n-1}$} \put(75,-10){$i_1$} \put(180,-55){$i_2$}
\put(126,-38){$g$} \put(170,-16){$\cup$}
\put(168,-30){$W^{\,n-1}$} \put(124,30){$\Psi_2$}
\put(235,0){$f|_U=\Psi_1\circ (i_1\times (\Psi_2\circ (g\times
i_2)))$} \put(145,10){\vector(0,1){60}}
\put(120,-64){\vector(0,1){60}} \put(175,-65){\vector(0,1){28}}
\put(50,-62){\vector(1,4){27}}
\put(90,80){$V^{p-1}\times_{\sigma_1} \Q^{\,n+3}_{\tilde c}$}
\put(88,66){$\cup$} \put(80,52){$W^{\,p-1}$}
\put(160,81){\vector(1,0){55}} \put(218,80){$\Q_c^{\,p+n+2}$}
\put(34,-80){$L^{p-1}_0\!\times_{\rho_1=\sigma_1\circ i_1}
\!\!((I\!\!\times_{\rho_3}\!J)\! \times_{\bar{\rho}=\sigma_2\circ
g}\!M^{n-1}_0)=U =L^p_0 \times_\rho M^n_0$} \put(175,86){$\Psi_1$}
\put(175,35){$\circlearrowright$} \put(225,-63){\vector(0,1){130}}
\end{picture}
\bigskip
\bigskip
\bigskip
\bigskip
\bigskip
\bigskip
\bigskip
\medskip
\end{itemize}

\end{theorem}

   In case $(iii)$ the \ii $g\colon\,I\times_{\rho_3}J\to W^{4}$
is neither a warped product $g=\Psi_3\circ(\alpha\times \beta)$,
where $\Psi_3\colon\,V^{1+k_1}\times_{\sigma_3}
\Q_{\hat{c}}^{\,1+k_2}\to \Q_{\tilde{c}}^{\,4}$ is a  warped
product representation with $k_1+k_2=2$ and $\alpha\colon\, I\to
V^{1+k_1}$ and $\beta\colon\, J\to \Q_{\hat{c}}^{\,1+k_2}$ are
unit speed curves with $\rho_3=\sigma_3\circ \alpha$, nor a
composition $H\circ G$ of such a warped product
$G=\Psi_3\circ(\alpha\times \beta)$, determined by a warped
product representation $\Psi_3\colon\,V^{1+k_1}\times_{\sigma_3}
\Q_{\hat{c}}^{\,1+k_2}\to \Q_{\tilde{c}}^{\,3}$ as before with
$k_1+k_2=1$, and an isometric immersion $H$ of an open subset
$W\supset G(I\times J)$ into $\Q_{\tilde{c}}^{\,4}$. It would be
interesting to exhibit an explicit example of such an isometric
immersion. Notice that it must satisfy the additional condition
$\sigma_2\circ g=(\rho_3\circ \pi_I)(\rho_2\circ \pi_J)$ for some
$\rho_2\in C^\infty(J)$.

    Cases $(i)-(iii)$ are disjoint. In fact, we will prove that
    under the assumptions of the theorem
    there are three distinct possible structures for the \sff of $f$, each of
    which corresponds to one of the cases in the statement.

   Notice that the conclusion of the theorem
remains unchanged under the apparently weaker assumption that the
subset of points of $N^{p+n}$ with constant sectional curvature
$c$ has empty interior.

  Theorem \ref{thm:main} does not hold without the assumption
that $n\geq 3$. In fact, we argue next that local \iis of the
round three-dimensional sphere $\Sf^{\,3}$ into $\R^5$ are
generically as in neither of the cases in the statement with
respect to any  local decomposition of $\Sf^{\,3}$ as a warped
product.

\begin{example}\label{esfera}
It was shown in  \cite{dt2} (cf.\  Corollary $4$ in \cite{dt3})
that local \iis of $\Sf^{\,3}$ into $\R^5$ that are {\it nowhere
compositions} (i.e., on no open subset they are compositions of
the umbilical inclusion into $\R^4$ with a local \ii of $\R^4$
into $\R^5$) are in correspondence with solutions $(V,h)$ on open
simply connected subsets $U_0\subset \R^3$ of the nonlinear system
of PDE's
$$
(I)\;\;\; \left\{  \begin{array}{l} (i)\; {\displaystyle\frac{\d
V_{ir}}{\d u_{j}} =
h_{ji}V_{jr}},\;\;\;\;\;(ii)\;{\displaystyle\frac{\d h_{ik}}{\d
u_j} = h_{ij}h_{jk}},
\vspace{2ex}\\
(iii)\;  {\displaystyle\frac{\d h_{ij}}{\d u_i} + \frac{\d
h_{ji}}{\d u_j} + \sum_{k}h_{ki}h_{kj} + V_{i3}V_{j3} =0},
\;\;\;\;i\neq j\neq k\neq i,\vspace{1.5ex} \\
\end{array} \right.
$$
called the generalized elliptic sinh-Gordon equation. Here
$V\colon\,U_0\to \Oes_1(3)$ is a smooth map taking values in the
group of orthogonal matrices with respect to the Lorentz metric of
signature $(+,+,-)$ and $x\in U_0\mapsto h(x)$ is a smooth map
such that $h(x)$ is an off-diagonal $(3\times 3)$-matrix for every
$x\in U_0$. More precisely, for any such \ii there exist a local
system of coordinates $(u_1,u_2,u_3)$, an orthonormal normal frame
$\{\xi_1,\xi_2\}$ and matrix functions $V$ and $h$ as above such
that \be\label{eq:ex1}
\!\!\!\!\!A_{\xi_r}X_i=V_{i3}^{-1}V_{ir}X_i,\;\,1 \leq  r\leq  2,
\;\,1\leq  i\leq 3,\ee and
\be\label{eq:ex1b} \nabla_{\d/\d
u_i}X_j=h_{ji}X_i, \;\,1 \leq i\neq j\leq 3, \ee where $X_i$ is a
unit vector field with $\d/\d u_i=V_{i3}X_i$. The compatibility
equations for $f$ are equivalent to system $(I)$. Conversely, any
solution $(V,h)$ of system $(I)$ on an open simply connected
subset $U_0\subset \R^3$ gives rise to such an \ii by means of the
fundamental theorem of submanifolds.

 By a theorem of Bourlet (see \cite{bo},
Th\'eor\`eme VIII), there exists one and only one analytic
solution $(V,h)$ of system $(I)$ in a neighborhood of an initial
value $u_0=(u_1^0,u_2^0,u_3^0)$ such that $V(u_0)\in \Oes_1(3)$
and such that $V_i^k$ and $h_{ij}$, $i<j$ (resp., $i>j$) reduce to
an arbitrarily given analytic function of $u_i$ (resp., $u_j$)
when the remaining variables take their initial values. Thus, for
a generic local analytic solution $(V,h)$ the functions $h_{ij}$
are nowhere vanishing; see the last section of \cite{dt3} for
explicit  \iis  with this property.

  It follows easily from  (\ref{eq:ex1}) that no such \ii admits a normal vector
  field whose shape operator has rank one. In particular,
  it can not be as in case $(ii)$.  Also, if $f$ is
as in case $(i)$ with respect to a decomposition
$U=L^{k_1}\times_{\rho} M^{k_2}$ of $U$ as a warped product, then
we must have that $k_2=1$. In fact,  otherwise  the second
fundamental form of $f$ would be given by
$$
\alpha(Y,Z)=\<Y_1,Z_1\>\eta_1+\<Y_2,Z_2\>\eta_2
$$
for some normal vector fields $\eta_1, \eta_2$ satisfying
$\<\eta_1,\eta_2\>=1=\|\eta_2\|$, where $Y_i, Z_i$, $1\leq i\leq
2$, are the components of $Y,Z$ according to the product
decomposition of $U$. This easily implies that the shape operator
with respect to a normal vector field orthogonal to $\eta_2$ has
rank one.  Thus, the distribution tangent to the second factor is
one-dimensional and invariant by all shape operators of $f$, and
hence it must be spanned by one of the vector fields $X_i$, $1\leq
i\leq 3$, say, $X_3$. In particular, this implies that the
distribution spanned by $X_1, X_2$ is totally geodesic, and hence
the functions $h_{31}$ and $h_{32}$ vanish everywhere by
$(\ref{eq:ex1b})$. Finally, we claim that the same holds if $f$ is
as in case $(iii)$. In effect, in this case $U$ splits as a warped
product
$$
U=L^1\times_{\rho}\Q_{\tilde{c}}^{\,2}
=L^1\times_\rho(J\times_{\bar{\rho}} M^1) =(L^1\times_{\rho}
J)\times_{(\rho\circ\pi_{L^1})\,\bar{\rho}}M^1,
$$
and $f$ is a warped product $f=g\times i$ with respect to the last
decomposition.  Thus, we have again that the one-dimensional
distribution tangent to $M^1$ is invariant by all shape operators
of $f$, and the same argument used in the preceding case proves
our claim. It follows that $f$ is  generically as in neither of
the cases in Theorem \ref{thm:main} with respect to any  local
decomposition of $\Sf^{\,3}$ as a warped product. \end{example}

 We now discuss some further results. The case of hypersurfaces is interesting in its own right. Although it
can be proved as a corollary of Theorem~\ref{thm:main}, it is
easier to derive it as an immediate consequence of Theorem
\ref{thm:no} and Proposition~\ref{cor:hyp} of Section $3$.

\begin{theorem}\label{thm:hyp} Assume that a warped product   $L^p\times_\rho
M^n$, $n\geq 2$, has no points with constant sectional curvature
$c$. Then any \ii $f\colon\,L^p\times_\rho M^n\to\Q_c^{\,p+n+1}$
is a warped product $f=\Psi\circ (F\times G)$, where $\Psi\colon\,
V^{p+k_1}\times_{\sigma}\Q_{\tilde{c}}^{\,n+k_2} \to
\Q_c^{\,p+n+1}$ is a warped product representation with
$k_1+k_2=1$ and $F\colon\, L^p\to V^{p+k_1}$, $G\colon\,
M^n\to\Q_{\tilde{c}}^{\,n+k_2}$ are isometric immersions.
\end{theorem}
Again, the preceding result is false if the assumption that $n\geq
2$ is dropped; rotation surfaces in $\R^3$ admit many isometric
deformations into nonrotational surfaces (cf.\ \cite{bi}).

  In deriving Theorem \ref{thm:main} we also obtain the
following result for the case of Riemannian products, which
extends Theorem $1$ in \cite{mo} in the case of products with two
factors. Therein, isometric immersions of Riemannian products with
arbitrarily many factors into Euclidean space were shown to split
as a product of isometric immersions under the assumptions that no
factor has an open subset of flat points and that the codimension
equals the number of factors. We point out that in the case of
Riemannian products the factors may change the roles. This
observation is applied several times throughout the paper.
\begin{theorem}\label{cor:prod2}
Let $f\colon\,L^p\times M^n\to\Q_c^{\,p+n+2}$ be an \ii of a
Riemannian  product. If $c=0$ assume that either $L^p$ or $M^n$
has dimension at least two and is free of flat points. If $c\neq
0$ assume that either $n\geq 3$ or $p\ge 3$. Then there exists an
open dense subset of $L^p\times M^n$ each of whose points lies in
an open product neighborhood $U=L^p_0\times M^n_0\subset L^p\times
M^n$ such that one of the following possibilities
holds:\vspace{.5ex}

\noindent Case  $c=0.$
\begin{itemize}
\item[$(i)$] There exist an orthogonal decomposition $\R^{\,p+n+2}=\R^{\,p+k_1}\times\R^{\,n+k_2}$
with $k_1+k_2=2$ and \iis $h_1\colon\,L^p_0\to\R^{\,p+k_1}$ and
\mbox{$h_2\colon\,M^n_0\to\R^{\,n+k_2}$} such that
\mbox{$f|_U=h_1\times h_2$}.

\vspace{3ex}
\begin{picture}(150,84)\hspace*{-10ex}
\put(100,30){$L^p_0\;\;\times \, M^n_0$} \put(112,55){$h_1$}
\put(160,32){\vector(2,3){30}} \put(150,55){$h_2$}
\put(178,50){$f|_U=h_1\times h_2$} \put(107,42){\vector(0,1){30}}
\put(145,42){\vector(0,1){30}} \put(102,80){$\R^{p+\!k_1}\!\times
\R^{n+k_2}=\R^{p+n+2}$}
\end{picture}
\vspace*{-5ex}

\item[$(ii)$] There exist an orthogonal decomposition
$\R^{\,p+n+1}=\R^{\,p+k_1}\times \R^{\,n+k_2}$, $k_1+k_2=1$, and
\iis $h_1\colon\,L^p_0\to \R^{\,p+k_1}$, $h_2\colon M^n_0\to
\R^{\,n+k_2}$ and $H\colon\,W\to \R^{\,p+n+2}$ of an open subset
$W\supset (h_1\times h_2)(U)$ of $\;\R^{\,p+n+1}$ such that
$f|_U=H\circ (h_1\times h_2)$.
\end{itemize}
\bigskip
\bigskip
\begin{picture}(150,84)\hspace*{-10ex}
\put(100,30){$L^p_0\;\;\times \, M^n_0$} \put(112,55){$h_1$}
\put(160,32){\vector(1,0){138}} \put(305,30){$\R^{p+n+2}$}
\put(143,55){$h_2$} \put(210,45){$\circlearrowright$}
\put(180,20){$f|_U=H\circ (h_1\times h_2)$}
\put(107,42){\vector(0,1){30}} \put(138,42){\vector(0,1){30}}
\put(190,80){\vector(3,-1){110}}
\put(102,80){$\R^{p+\!k_1}\!\times \R^{n+k_2} \!\supset\! W$}
\put(242,65){$H$}
\end{picture}
\vspace*{-3ex}

\noindent Case $c\neq 0.$\vspace{.5ex}
\begin{itemize}
\item[$(i)$]
There exist an  embedding $\Phi\colon\,\Q_{c_1}^{\,p+k_1}\times
\Q_{c_2}^{\,n+k_2} \to \Q_c^{\,p+n+2}$ as an extrinsic Riemannian
product with $k_1+k_2=1$, and \iis $h_1\colon\,L^p_0\to
\Q_{c_1}^{\,p+k_1}$ and \mbox{$h_2\colon\,M^n_0\to
\Q_{c_2}^{\,n+k_2}$} such that $f|_U=\Phi\circ(h_1\times h_2)$.

\bigskip
\bigskip
\begin{picture}(150,84)\hspace*{-15ex}
\put(100,30){$L^p_0\;\;\times \; M^n_0$} \put(112,55){$h_1$}
\put(160,35){\vector(2,1){85}} \put(145,55){$h_2$}
\put(205,45){$f|_U=\Phi\circ(h_1\times h_2)$}
\put(107,42){\vector(0,1){30}} \put(140,42){\vector(0,1){30}}
\put(100,80){$\Q_{c_1}^{\,p+k_1}\!\times\! \Q_{c_2}^{\,n+k_2}$}
\put(167,83){\vector(1,0){76}} \put(180,60){$\circlearrowright$}
\put(200,87){$\Phi$} \put(250,80){$\Q_c^{\,p+n+2}$}
\end{picture}
\vspace*{-5ex}

\item[$(ii)$] There exist an  embedding $\Phi\colon\,\Q_{c_1}^{\,p}\times \Q_{c_2}^{\,n}\to \Q_c^{\,p+n+1}$ as an extrinsic Riemannian product, local isometries $i_1\colon\, L_0^p\to \Q_{c_1}^{\,p}$ and
$i_2\colon\, M_0^n\to \Q_{c_2}^{\,n}$,  and an \ii $H\colon\,W\to
\Q_c^{\,p+n+2}$ of an open subset $W\supset \Phi\circ (i_1\times
i_2)(U)$ of  $\Q_c^{\,p+n+1}$ such that $f|_U=H\circ \Phi\circ
(i_1\times i_2)$.

\bigskip
\bigskip
\bigskip

\begin{picture}(150,84)\hspace*{-18ex}
\put(100,30){$L_0^p\times M_0^n$} \put(110,55){$i_1$}
\put(155,32){\vector(1,0){135}} \put(295,30){$\Q_c^{\,p+n+2}$}
\put(142,55){$i_2$} \put(170,20){$f|_U=H\circ\Phi\circ(i_1\times
i_2)$} \put(107,42){\vector(0,1){30}}
\put(137,42){\vector(0,1){30}} \put(100,80){$\Q_{c_1}^{\,p}\times
\Q_{c_2}^{\,n}$} \put(155,83){\vector(1,0){135}}
\put(295,80){$W\subset\Q_c^{\,p+n+1}$} \put(302,55){$H$}
\put(298,72){\vector(0,-1){30}} \put(215,90){$\Phi$}
\put(215,55){$\circlearrowright$}
\end{picture}
\vspace*{-3ex}
\end{itemize}

\end{theorem}

        As an example showing that for $c\neq 0$ the
assumption  that either $n\geq 3$ or $p\ge 3$ is indeed necessary,
we may take any local \ii of $\R^3$ into $\Sf^{\,5}$ that is not a
product $ \alpha\times g\colon\, I\times V\to \Sf^{\,1}(r_1)\times
\Sf^{\,3}(r_2),\;\; r_1^2+r_2^2=1,$ where $\alpha\colon\,I\to
\R^2$ is a unit speed parametrization of an open subset of a
circle of radius~$r_1$ and \mbox{$g\colon\,V\to  \Sf^{\,3}(r_2)$}
is an isometric immersion of an open subset $V\subset\R^2$. Recall
that local \iis of $\R^3$ into $\Sf^{\,5}$ were shown in
\cite{ten} to be in correspondence with solutions on simply
connected open subsets of $\R^3$ of the so-called generalized wave
equation. As in the previous discussion on local \iis of
$\Sf^{\,3}$ into $\R^5$, one may easily argue that the class of
local \iis of $\R^3$ into $\Sf^{\,5}$ that are given as products
as just described is only a rather special subclass of the whole
class of such isometric immersions.\vspace{1ex}

   We now give precise statements of the applications of
Theorem \ref{thm:main} referred to at the end of the introduction.
Recall that an \ii $F\colon\,L^p\to\Q_c^{\,p+m}$ is said to be
{\em locally rigid\/} if it is rigid restricted to any open subset
of $L^p$.

\begin{corollary} \label{cor:st1}
Let $L^p$ be a Riemannian manifold no open subset of which can be
isometrically immersed in $\Q_c^{\,p+1}$. If
$f\colon\,L^p\times_{\rho}M^n\to \Q_c^{\,p+n+2}$, $n\geq 3$, is an
isometric immersion, then there exist a warped product
representation  \mbox{$\Psi\colon V^{p+2}\!\times\!
_{\sigma}\Q_{\tilde{c}}^{\,n}\to \Q_c^{\,p+n+2}$}, an \ii
$F\colon\,L^p\to V^{p+2}$ and a local isometry \mbox{$i\colon
M^n\to\Q_{\tilde{c}}^{\,n}$} such that $f=\Psi\circ (F\times i)$.
In particular, if $L^p$ is a Riemannian manifold that admits a
locally rigid \ii $F\colon\,L^p\to \Q_c^{\,p+2}$, then the
preceding conclusion holds and, in addition, the \ii $f$ is also
locally rigid.
\end{corollary}
\proof By Theorem \ref{thm:main}, any \ii
$f\colon\,L^p\times_{\rho}M^n\to \Q_c^{\,p+n+2}$, $n\geq 3$, must
be locally as in one of the three cases described in its
statement. However, under the assumption that $L^p$ has
 no open subset that can be isometrically
immersed in $\Q_c^{\,p+1}$, it follows that $f|_U$ can not be as
in case $(ii)$ on any open subset $U=L^p_0\times M^n_0\subset
L^p\times M^n$, for there can  not exist by that assumption any
isometric immersion $h_1\colon\,L_0^p\to \Q_c^{\,p+k_1}$ with
$0\leq k_1\leq 1$. Moreover, $f|_U$ can not be as in case $(iii)$
either on any such open subset, for in that case $L^p_0$ would
have constant sectional curvature $c$, and hence it would admit
locally an \ii into $\Q_c^{\,p+1}$. Therefore $f$ must be globally
as in case~$(i)$. The last assertion is now clear.
\vspace{1ex}\qed

    We say that a Riemannian manifold can be locally isometrically immersed in $\Q_c^{\,\ell}$ if each point
has an open neighborhood that admits an isometric immersion into
$\Q_c^{\,\ell}$. Arguing in a similar way as in the proof of
Corollary \ref{cor:st1} yields the following.

\begin{corollary} Let $L^p$ be a Riemannian manifold
that cannot be locally isometrically immersed in $\Q_c^{\,p+2}$.
Then $L^p\times_{\rho}M^n$ can not be locally isometrically
immersed in $\Q_c^{\,p+n+2}$ for any Riemannian manifold $M^n$ of
dimension $n\geq 3$ and any warping function $\rho$.
\end{corollary}

  In view of N\"olker's result it is also natural to study
\iis of warped products into space forms in the light of the
following definition.\vspace{1,5ex}

\noindent{\bf Definition.} Let $f\colon\,L^p\times _{\rho}M^n\to
\Q_c^{\,\ell}$ be an \ii of a warped product. Given  $z\in
L^p\times_{\rho}M^n$, we denote
$$
\a({\mathcal H}_z,{\mathcal V}_z) =\spa \{\a(Y,V):Y\in\H_z\an V\in
\V_z\}.
$$
We say that  the immersion $f$ at $z$ is of {\it type}
\begin{itemize}
\item[$(A)$] if $\dim \a({\mathcal H}_z, {\mathcal V}_z)=0$,
i.e., $\a(Y,V)=0$ for all $Y\in {\mathcal H}_z$ and $V\in
{\mathcal V}_z$,
\item[$(B)$] if $\dim \a({\mathcal H}_z,{\mathcal V}_z)=1$,
\item[$(C)$] if $\dim \a({\mathcal H}_z,{\mathcal V}_z)
\geq 2$.\end{itemize}

    Notice that type $A$ is the case of N\"olker's
decomposition theorem.  Therefore, it is a natural problem to
determine the \iis that are everywhere of type~$B$. In the
following section we obtain a complete solution to this problem in
the codimension two case under the assumption that $n\geq 3$ (for
the case of Riemannian products it is enough to assume that
$p+n\geq 3$); see the two paragraphs before
Proposition~\ref{pro:cij0}. The proof of Theorem~\ref{thm:main} is
then accomplished as follows: type $C$ is excluded by
Proposition~\ref{casec}, type $A$ corresponding to Theorem
\ref{thm:no} gives case $(i)$, and type $B$ splits into two
subcases $B_1$ and $B_2$ handled in Propositions \ref{pro:cij0}
and \ref{le:B_5}, respectively,  which correspond to the cases
$(ii)$ and $(iii)$. Similarly for the proof of Theorem
\ref{cor:prod2}: type $C$ is excluded by
Corollary~\ref{casecprod}, type $A$ corresponding to Theorem
\ref{thm:molzan} gives subcase $(i)$ in both cases $c=0$ and
$c\neq 0$, type $B_1$ handled in Corollary \ref{cor:rpb} gives
subcase $(ii)$ in either case, and type $B_2$ is excluded in
either case by Corollary~\ref{escolio} and Corollary
\ref{cor:prodb2}, respectively.

\section{Immersions of type $B$}

 Our main goal in this section is to provide a complete local classification of \iis  $f\colon\, L^p\times_\rho M^n\to\Q_c^{\,p+n+2}$  that are everywhere of type $B$
under the assumption that $n\geq 3$.  A similar classification for
the special case of \iis of Riemannian products is also given, for
which it is enough to assume $p+n\ge 3$.\vspace{1ex}

First we determine the pointwise structure of the second
fundamental forms of \iis of type $B$,  starting with some general
facts that are valid in arbitrary codimension.

\begin{lemma}\label{prop:rankone} Let
$f\colon\,L^p\times_\rho M^n\to \Q_c^{\,\ell}$ be an \ii of a
warped product. Assume that $f$ is not of type $A$ at a point
$z\in N^{p+n}=L^p\times_\rho M^n$ and that for every $Y\in
{\mathcal H}_z$ the linear map \be\label{eq:bx} B_Y\colon\,
{\mathcal V}_z \to T_z^\perp N, \;\;\;V\mapsto\a(Y,V), \ee
satisfies  $\rank B_Y\leq 1$. Then there exists a unit vector
$e\in {\mathcal V}_z$, uniquely determined up to its sign, such
that \be\label{eq:r1a} \a(Y,V)=\<V,e\>\a(Y,e) \ee and
\be\label{eq:r1b} \a(V,W)-\<V,e\>\<W,e\>\a(e,e) \perp\a({\mathcal
H}_z,{\mathcal V}_z) \ee for all $Y\in {\mathcal H}_z$ and $V,W\in
{\mathcal V}_z$.
\end{lemma}

\proof Let $X\in {\mathcal H}_z$ be such that $\rank B_X = 1$.
Then ${\mathcal D}(X):=\ker B_X$ has codimension $1$ in ${\mathcal
V}_z$. Let $e\in {\mathcal V}_z$ be one of the unit vectors
perpendicular to ${\mathcal D}(X)$ and write $B_X e=\lambda \xi$
where $\lambda\neq 0$ and $\xi\in T_z^\perp N$ is a unit vector.
Let $Y\in {\mathcal H}_z$ and $V\in {\mathcal D}(X)$ be arbitrary
vectors. Then (\ref{eq:ge33}) implies
\be\label{eq:r8}\begin{array}{l} \<B_Xe,B_YV\>=\<\a(X,e),
\a(Y,V)\> =\<\a(X,V), \a(Y,e)\>\vspace{1ex}\\\hspace{13.5ex}
=\<B_X V,\a(Y,e)\>=0. \end{array}\ee Now consider the linear map
$B_{X+tY}$ for arbitrary $t\in \R$. By assumption its rank is at
most $1$. Therefore the vectors $B_{X+tY} e=\lambda\xi +tB_Y e$
and $B_{X+tY} V=tB_Y V$ are linearly dependent, and hence
$$
\<B_{X+tY}e,B_{X+tY}e\>\<B_{X+tY}V,B_{X+tY} V\>
-\<B_{X+tY}e,B_{X+tY}V\>^2=0.
$$
As the left hand side of this equation is a polynomial
$\sum_{i=2}^4 a_it^i$, its coefficients must vanish; in
particular, because of (\ref{eq:r8}) we obtain
$0=a_2=\lambda^2\|B_Y V\|^2$. Hence $B_Y|_{{\mathcal D}(X)}=0$,
and (\ref{eq:r1a}) follows.\vspace{1ex}

  By means of (\ref{eq:ge44}) we derive
$$
\<\a(V,W),\a(Y,e)\>=\<\a(Y,V),\a(e,W)\>
=\<V,e\>\<\a(e,W),\a(Y,e)\>.
$$
Applying this result to $\a(W,e)$ instead of $\a(V,W)$ we obtain
$$
\<\a(V,W)-\<V,e\>\<W,e\>\a(e,e),\a(Y,e)\>=0,
$$
which implies (\ref{eq:r1b}) in view of (\ref{eq:r1a}). \qed

\begin{lemma}\label{le:typeb} Let
$f\colon\,L^p\times_\rho M^n\to\Q_c^{\,\ell}$ be an \ii of a
warped product. Assume that $f$ is of type $B$ at $z\in
N^{p+n}=L^p\times_\rho M^n$. Then there exist unique, up to their
signs, unit vectors $X\in {\mathcal H}_z$, $e\in {\mathcal V}_z$
and $\xi\in T_z^\perp N$,  and $\beta,\lambda, \gamma\in\R$ with
$\lambda\neq 0$ such that
\begin{eqnarray}
\<\a(Y,Z),\xi\>\!\!\!&=&\!\!\!\beta \<Y,X\>\<Z,X\>,
\label{eq:beta}\\
\a(Y,V)\!\!\!&=&\!\!\!\lambda\<Y,X\>\<V,e\>\xi,
\label{eq:lambda}\\
\<\a(V,W),\xi\>\!\!\!&=&\!\!\!\gamma \<V,e\>\<W,e\>,
\label{eq:gamma}\end{eqnarray}
\be\label{eq:tildep}\begin{array}{l} \<\tilde{P}\a(Y,Z),
\tilde{P}\a(V,W)
-\<V,W\>\tilde{P}\a(e,e)\>\vspace{1ex}\\\hspace{25ex}=(\beta\gamma-\lambda^2)
\<Y,X\>\<Z,X\>\<PV,PW\>,\end{array}\ee where
$\tilde{P}\colon\,T_z^\perp N\to T_z^\perp N$ and
$P\colon\,{\mathcal V}_z\to{\mathcal V}_z$ denote the orthogonal
projections  onto the subspaces $\{\xi\}^\perp\subset T_z^\perp N$
and $\{e\}^\perp\subset {\mathcal V}_z$, respectively. Moreover,
if $N^{p+n}=L^p\times M^n$ is a Riemannian product then \be
\label{eq:tildep2}\begin{array}{l}
\<\tilde{P}\a(Y,Z),\tilde{P}\a(V,W)\>+ (\beta\gamma
-\lambda^2)\<Y,X\>\<Z,X\>\<V,e\>\<W,e\>\vspace{1ex}\\\hspace{40ex}
+ c\<Y,Z\>\<V,W\>=0.\end{array} \ee
\end{lemma}

\proof Let $\xi\in T_z^\perp N$ be a unit vector such that
$\a({\mathcal H}_z,{\mathcal V}_z)=\R\,\xi$. Given $Y\in {\mathcal
H}_z$, then $B_Y$ takes its values in $\R\, \xi$, and hence $\rank
B_Y\leq 1$. Thus, we may apply Lemma~\ref{prop:rankone}.  On the
other hand, since the linear map ${\mathcal H}_z\to T_z^\perp N$
defined by $Y\mapsto \a(Y,e)$ also takes its values in $\R\, \xi$,
it follows that, up to sign, there exists exactly one unit vector
$X\in {\mathcal H}_z$ perpendicular to its kernel. Set
$\gamma=\<\alpha(e,e),\xi\>$, $\lambda=\<\a(X,e),\xi\>$ and
$\beta=\<\alpha(X,X),\xi\>$. Notice that $\lambda\neq 0$ because
$f$ is of type $B$ at $z$. We obtain (\ref{eq:gamma})  from
(\ref{eq:r1b}), whereas (\ref{eq:lambda}) follows from
(\ref{eq:r1a}) and $\a(Y,e)=\<Y,X\>\a(X,e)=\lambda\<Y,X\>\xi$.
Using (\ref{eq:ge22}) we obtain
$$
\lambda\<\a(Y,Z),\xi\>= \< \a(Y,Z),\a(X,e)\>=\< \a(X,Z),\a(Y,e)\>
=\lambda\<Y,X\>\<\a(X,Z),\xi\>,
$$
and applying this result to $\a(Z,X)$ instead of $\a(Y,Z)$ we end
up with (\ref{eq:beta}).

 We obtain from (\ref{eq:ge11}),
(\ref{eq:beta}),  (\ref{eq:lambda}) and (\ref{eq:gamma}) that
$$
\begin{array}{l}
\<V,W\>\<\nabla_Y\eta-\<Y,\eta\>\eta-cY,Z\>\vspace{1ex}\\
\hspace*{12ex}=\<\a(Y,Z),\a(V,W)\>
-\<\a(Y,W),\a(Z,V)\>\vspace{1ex}\\
\hspace*{12ex}=(\beta\gamma
-\lambda^2)\<Y,X\>\<Z,X\>\<V,e\>\<W,e\>
+\<\tilde{P}\a(Y,Z),\tilde{P}\a(V,W)\>.
\end{array}
$$
This yields (\ref{eq:tildep2}) if $N^{p+n}=L^p\times M^n$ is a
Riemannian product. In the general case,   putting $W=V=e$ we get
$$
\<\nabla_Y\eta-\<Y,\eta\>\eta-cY,Z\>
=(\beta\gamma-\lambda^2)\<Y,X\>\<Z,X\>
+\<\tilde{P}\a(Y,Z),\tilde{P}\a(e,e)\>.
$$
The two preceding equations yield
$$
\begin{array}{l}
\<\tilde{P}\a(Y,Z), \tilde{P}\a(V,W)
-\<V,W\>\tilde{P}\a(e,e)\>\vspace{1ex}\\\hspace*{12ex}=(\beta\gamma
-\lambda^2)\<Y,X\>\<Z,X\>(\<V,W\>-\<V,e\>\<W,e\>),
\end{array}
$$
which coincides with (\ref{eq:tildep}). \qed

\begin{remark and definition}\label{re:typeb}  Equations (\ref{eq:beta}), (\ref{eq:lambda})
and (\ref{eq:gamma}) are equivalent to \be\label{eq:shape1} A_\xi
Y=\<Y,X\>(\beta X+\lambda e),\;\;\;\; A_\xi V=\<V,e\>(\lambda
X+\gamma e), \ee and \be\label{eq:shape11} \tilde{P}\a(Y,V)=0. \ee
In particular, it follows from (\ref{eq:shape1}) that the rank of
$A_\xi$ at $z$ is either $1$ or $2$, according to
$\beta\gamma-\lambda^2$ being zero or not. We say accordingly that
$f$ is of {\em type\/}~$B_1$ or of {\em type\/} $B_2$ at $z$.
\end{remark and definition}

 We now show that in the case of hypersurfaces
 $f\colon\,L^p\times_\rho M^n\to\Q_c^{\,p+n+1}$, $n\geq 2$,
 only types $A$ and $B_1$ can occur pointwise.

\begin{proposition}\label{cor:hyp}
Let $f\colon\,L^p\times_\rho M^n\to\Q_c^{\,p+n+1}$, $n\geq 2$, be
an \ii of a warped product. Then, at any point $z\in
N^{p+n}=L^p\times_\rho M^n$ either $f$ is  of type $A$ or of type
$B_1$. Moreover, in the latter case  $N^{p+n}$ has constant
sectional curvature $c$ at $z$.
\end{proposition}

\proof Assume that $f$ is not of type $A$ at $z$. Since $n\geq 2$,
we may choose a unit vector $V\in \{e\}^\perp\subset {\mathcal
V}_z$. Applying (\ref{eq:tildep}) for $W=V$ and $Z=Y=X$, and using
that $\tilde{P}=0$, it follows that $\beta\gamma-\lambda^2=0$.
Therefore $f$ is of type $B_1$ at $z$. The last assertion follows
from the Gauss equation of $f$.\qed \vspace{1ex}

   Theorem \ref{thm:hyp} now follows by putting together the preceding result and
    Theorem \ref{thm:no}. For Riemannian products we obtain the following corollary.

\begin{corollary}\label{cor:prod21}
Let $f\colon\,L^p\times M^n\to\Q_c^{\,p+n+1}$ be an \ii of a
Riemannian product. Assume that $p+n\ge 3$ and, if $c=0$, that
either $L^p$ or $M^n$, say, the latter, has dimension at least two
and is free of flat points. Then $f$ is of type $A$ everywhere and
we have:
\begin{itemize}
\item[$(i)$] If $c=0$ there exist an orthogonal decomposition $\R^{p+n+1}=\R^{p}\oplus \R^{n+1}$, a local isometry $i\colon\,L^p\to \R^p$ and an \ii
$h\colon\,M^n\to \R^{n+1}$ such that \mbox{$f=i\times h$}.
\item[$(ii)$] If $c\neq 0$  there exist an embedding
$\Phi\colon\,\Q_{c_1}^{\,p}\times \Q_{c_2}^{\,n}\to
\Q_c^{\,p+n+1}$ as an extrinsic Riemannian product and local
isometries $i_1\colon\, L^p\to \Q_{c_1}^{\,p}$ and $i_2\colon\,
M^n\to \Q_{c_2}^{\,n}$ such that $f = \Phi\circ (i_1\times i_2)$.
\end{itemize}
\end{corollary}

 From now on we consider  \iis
 $f\colon\,L^p\times_\rho M^n\to\Q_c^{\,p+n+2}$. Assume that $f$ is of
type $B$ at a point $z\in N^{p+n}=L^p\times_\rho M^n$ and let
$X,e,\xi,\beta,\lambda$ and $\gamma$ be as in Lemma {\em
\ref{le:typeb}}. Choose one of  the unit vectors $\tilde{\xi}\in
T_z^\perp N$  perpendicular to $\xi$, and define the symmetric
bilinear forms
$$
\begin{array}{c}
\tilde{\beta}\colon\,{\mathcal H}_z \times {\mathcal H}_z\to
\R,\;\;\;(Y,Z)\mapsto
\<\a(Y,Z),\tilde{\xi}\>\vspace{2ex},\\
\;\;\;\tilde{\gamma}\colon\,{\mathcal V}_z \times {\mathcal
V}_z\to \R,\;\;\;\;(V,W) \mapsto \<\a(V,W),\tilde{\xi}\>.
\end{array}
$$
Set also $\tilde{\beta}_0=\tilde{\beta}(X,X),\;
\tilde{\gamma}_0=\tilde{\gamma}(e,e)$ and
$\,\tilde{\delta}_0:=\tilde{\beta}_0
\tilde{\gamma}_0+\beta\gamma-\lambda^2$. Then Lemma~\ref{le:typeb}
can be strengthened as follows.

\begin{proposition}\label{cor:typeb} Let $f\colon\,L^p\times_\rho M^n\to\Q_c^{\,p+n+2}$
be an isometric immersion. Assume that $f$ is of type $B$ at $z\in
N^{p+n}=L^p\times_\rho M^n$. With the preceding notations we have:
\begin{itemize}
\item[$(i)$] If $f$ is of type $B_1$ at $z$, then
one of the following (not exclusive) possibilities holds:
$$
\tilde{\gamma}(V,W)=\<V,W\>\tilde{\gamma}_0
\;\;\;\;\;\;\mbox{or}\;\;\;\;\;  \tilde{\beta}= 0.
$$
\item[$(ii)$] If $n\geq 2$ and $f$ is of type $B_2$ at $z$, then
\be\label{eq:B_1} \tilde{\beta}(Y,Z)
=\<Y,X\>\<Z,X\>\tilde{\beta}_0
\;\;\;\mbox{with}\;\;\;\tilde{\beta}_0\neq 0,\; \mbox{and} \ee
\be\label{eq:B_2} \tilde{\beta}_0\,\tilde{\gamma}(V,W)
=\<V,W\>\tilde{\delta}_0 -(\beta\gamma -\lambda^2)\<V,e\>\<W,e\>.
\ee
\end{itemize}
\end{proposition}

\proof Equation (\ref{eq:tildep}) now reads
\be\label{eq:tildebeta} \tilde{\beta}(Y,Z)(\tilde{\gamma}(V,W)
-\<V,W\>\tilde{\gamma}_0)
=(\beta\gamma-\lambda^2)\<Y,X\>\<Z,X\>\<PV,PW\>. \ee If
$\beta\gamma-\lambda^2=0$ then the preceding equation proves
assertion $(i)$.   Choosing $Y=Z=X$ in (\ref{eq:tildebeta}) yields
\be\label{eq:xx} \tilde{\beta}_0\,(\tilde{\gamma}(V,W)
-\<V,W\>\tilde{\gamma}_0) =(\beta\gamma-\lambda^2)\<PV,PW\>. \ee
Using (\ref{eq:xx}), $n\geq 2$  and $\beta\gamma-\lambda^2\neq 0$
we derive $\tilde{\beta}_0\neq 0$ and from (\ref{eq:tildebeta}) it
follows that $\tilde{\beta}(Y,Z)=\<Y,X\>\<Z,X\>\tilde{\beta}_0$.
Finally, (\ref{eq:xx})  also yields (\ref{eq:B_2}).
\qed\vspace{1.5ex}

Taking into account (\ref{eq:tildep2}) we have the following
additional information in the case of Riemannian products.

\begin{corollary}\label{cor:typebprod}
Let $f\colon\,L^p\times M^n\to\Q_c^{\,p+n+2}$ be an \ii of a
Riemannian product. Assume that $f$ is of type $B$ at $z\in
L^p\times M^n$. Then, with the preceding notations, we have that
$\tilde{\delta}_0=-c$ and, in addition:
\begin{itemize}
\item[$(i)$] If $f$ is of type $B_1$ at $z$, then
$$
\tilde{\beta}(Y,Z)=\<Y,Z\>\tilde{\beta}_0,\;\;\;\;
\tilde{\gamma}(V,W)=\<V,W\>\tilde{\gamma}_0\;\;\;\;\;
\mbox{and}\;\;\;\;\;\tilde{\beta}_0\tilde{\gamma}_0+c=0.
$$

\item[$(ii)$] If $p,n\geq 2$ and $f$ is of type $B_2$ at $z$, then $c=0$,
 \be\label{eq:B_2prod} \tilde{\beta}(Y,Z)=\<Y,X\>\<Z,X\>\tilde{\beta}_0,\;\;\;
 \mbox{and}\;\;\;\;
\tilde{\gamma}(V,W)=\<V,e\>\<W,e\>\tilde{\gamma}_0.
 \ee
\end{itemize}
\end{corollary}

   In the remaining of this section we make a detailed
study of \iis $f\colon \,L^p\times_\rho M^n\to\Q_c^{\,p+n+2}$,
$n\geq 2$, that are everywhere of type $B$.  In this case we may
choose smooth unit vector fields $X$, $e$ and $\xi$ (and hence a
smooth unit  normal vector field $\tilde{\xi}$ orthogonal to
$\xi$), and smooth functions $\beta,\lambda$ and $\gamma$ that
satisfy pointwise the conditions of Lemma \ref{le:typeb}. In
Propositions \ref{pro:cij0} and \ref{le:B_5} below we classify
\iis of types $B_1$ and $B_2$, respectively, the latter only for
$n\geq 3$. In Corollary \ref{le:b6} we determine the special
subclass of \iis of type $B_2$ for which $\tilde{\delta}_0$  in
(\ref{eq:B_2}) is everywhere vanishing. Isometric immersions of
type $B_1$ of Riemannian products  with dimension $p+n\geq 3$ are
classified in Corollary \ref{cor:rpb}. In Corollary
\ref{cor:prodb2} we show that there exists no \ii $f\colon
\,L^p\times M^n\to\Q_c^{\,p+n+2}$, $p+n\geq 3$, of type $B_2$ if
$c\neq 0$ and in Corollary \ref{escolio} we classify such \iis for
$c=0$. This yields a local classification of \iis of type $B$ in
codimension two of warped products $L^p\times_{\rho} M^n$ for
which $n\geq 3$, as well as of Riemannian products $L^p\times M^n$
for which $p+n\geq 3$.

\begin{proposition}\label{pro:cij0} Let
$f\colon\,N^{p+n}=L^p\times_{\rho} M^n\to \Q_c^{\,p+n+2}$, $n\ge
2$,
 be an isometric embedding of type $B_1$. Then $f$ is a composition $H\circ g$
 of isometric immersions, where \mbox{$g=\Psi\circ(h_1\times h_2)$} is a warped~product of
\iis  determined by a warped product representation
$\Psi\colon\,V^{p+k_1}\times_{\sigma} \Q_{\tilde{c}}^{\,n+k_2} \to
\Q_c^{\,p+n+1}$, $k_1+k_2=1$, and $H\colon\, W\to \Q_c^{\,p+n+2}$
is an \ii of an open subset~$W\supset g(N^{p+n})$
of~$\,\Q_c^{\,p+n+1}$ $($see the diagram in
Theorem~\ref{thm:main}-$(ii))$.
\end{proposition}

\proof   Set  $E_0=(\beta X+\lambda e)/\|\beta X+\lambda e\|$.
Since $\beta\gamma-\lambda^2=0$ we have that $\lambda(\beta
X+\lambda e)=\beta(\lambda X+\gamma e)$, and hence
(\ref{eq:shape1}) yields \be\label{eq:axi}
A_{\xi}E=(\beta+\gamma)\<E,E_0\>E_0\fall E\in\Gamma(TN). \ee
Observe that $\beta+\gamma\neq 0$ because
$\beta\gamma-\lambda^2>0$.  The Gauss equation for $f$ and the
fact that $A_\xi$ has rank $1$ imply that $A_{\tilde{\xi}}$
satisfies the Gauss equation for an isometric immersion of
$N^{p+n}$ into $\Q_c^{\,p+n+1}$. We claim that it also satisfies
the Codazzi equation for such an isometric immersion. Define a
connection one-form $\omega$ on $N^{p+n}$ by
$\omega(E)=\<\nabla^\perp_E\xi,\tilde{\xi}\>$. By the Codazzi
equation for $f$ we have \be\label{eq:codazzi}
(\nabla_{E_1}A_{\tilde{\xi}})E_2 -(\nabla_{E_2}A_{\tilde{\xi}})E_1
=\omega(E_2)A_{\xi}E_1-\omega(E_1)A_{\xi}E_2. \ee The following
fact and (\ref{eq:axi}) imply that the right hand side of
(\ref{eq:codazzi}) vanishes, and the claim follows.

\begin{fact}\label{claim} {\em The one-form $\o$ satisfies $\omega(E)=\<E,E_0\>\omega(E_0)$ for all $E\in\Gamma(TN)$, or equivalently, $\omega(E)=0$
for all $E\in\Gamma(\ker A_{\xi})$. }\end{fact}

    In proving Fact \ref{claim} it is useful to observe that

     \be\label{eq:eqn1}
\ker A_{\xi(z)}=\{E_0(z)\}^\perp =\spa\{\lambda\<V,e\>Y
-\beta\<Y,X\>V\colon\,Y\in {\mathcal H}_z, V\in {\mathcal V}_z\}
\ee \be\label{eq:eqn2} \hspace*{22.2ex}=\spa\{\lambda\<Y,X\>V
-\gamma\<V,e\>Y\colon\,Y\in {\mathcal H}_z, V\in {\mathcal V}_z\}.
\ee

   By Proposition \ref{cor:typeb}-$(i)$, at each point
$z\in N^{p+n}$ either $A_{\tilde{\xi}}|_{{\mathcal V}_z}=
\tilde{\gamma}_0\,id$, where $id$ denotes the identity tensor, or
$A_{\tilde{\xi}}|_{{\mathcal H}_z}= 0$. Since $\o(E)$ is a
continuous function, it suffices to prove that $\o(E)(z)=0$ at
points $z\in N^{p+n}$ that are contained in an entire neighborhood
$U\subset N^{p+n}$ in which one of the preceding possibilities
holds everywhere. \vspace{1ex}

\noindent{\it Case $A_{\tilde{\xi}}|_{\mathcal V}=
\tilde{\gamma}_0\,id$.\/} For $Y\in {\mathcal L}(L^p)$ and $V,W\in
{\mathcal L}(M^n)$ we obtain using (\ref{eq:c2}) and~(\ref{eq:c4})
that
\begin{eqnarray}
\<(\nabla_Y A_{\tilde{\xi}})V -(\nabla_V A_{\tilde{\xi}})Y,W\>
\!\!\!&
=&\!\!\!\<\nabla_Y A_{\tilde{\xi}}V-\nabla_V A_{\tilde{\xi}}Y-A_{\tilde{\xi}}[Y,V],W\>\\
\!\!\!&=&\!\!\!\<\nabla_Y(\tilde{\gamma}_0V),W\>+\<A_{\tilde{\xi}}Y,
\nabla_V W\>\\
\!\!\!&=&\!\!\! (Y(\tilde{\gamma}_0) +\<(A_{\tilde{\xi}}
-\tilde{\gamma}_0\,id)Y,\eta\>)\<V,W\>.
\end{eqnarray}
On the other hand, by the Codazzi equation we have
\begin{eqnarray*}
\<(\nabla_Y A_{\tilde{\xi}})V -(\nabla_V A_{\tilde{\xi}})Y,W\>
\!\!\!&=&\!\!\!\<\omega(V)A_\xi Y-\omega(Y)A_\xi V, W\>\\
\!\!\!&=&\!\!\!\<W,e\>\omega(\lambda\<Y,X\>V -\gamma\<V,e\>Y),
\end{eqnarray*}
where in the second equality we have used (\ref{eq:shape1}). Thus
\be\label{eq:claim} (Y(\tilde{\gamma}_0)+\<(A_{\tilde{\xi}}
-\tilde{\gamma}_0\,id)Y,\eta\>)\<V,W\>
=\<W,e\>\omega(\lambda\<Y,X\>V-\gamma\<V,e\>Y) \ee for all $Y\in
{\mathcal L}(L^p)$ and $V,W\in {\mathcal L}(M^n)$. As these
equations are tensorial, they are also valid for arbitrary
horizontal (resp., vertical) vector fields $Y$ (resp., $V, W$).
In particular, if we apply (\ref{eq:claim}) for $W=V$ orthogonal
to $e$ we obtain that the expression between parenthesis on the
left-hand-side vanishes. Then, for $W=e$ this yields
$$
\omega(\lambda\<Y,X\>V-\gamma\<V,e\>Y)=0,
$$
thus proving Fact \ref{claim} by (\ref{eq:eqn2}) in this case.
\vspace{1ex}

\noindent{\it Case $A_{\tilde{\xi}}|_{\mathcal H}= 0$.\/}
 Using (\ref{eq:c1})
we obtain for $Y,Z\in {\mathcal L}(L^p)$ and $V\in{\mathcal
L}(M^n)$ that
\begin{eqnarray*}\<(\nabla_Y A_{\tilde{\xi}})V,Z\>
\!\!\!&=&\!\!\!\<\nabla_Y A_{\tilde{\xi}} V,Z\>
-\<A_{\tilde{\xi}}\nabla_Y V,Z\>\\
\!\!\!&=&\!\!\! Y\<A_{\tilde{\xi}} V,Z\> -\<A_{\tilde{\xi}}
V,\nabla_Y Z\> -\<\nabla_Y V, A_{\tilde{\xi}}Z\>=0
\end{eqnarray*}
and, analogously, that
$$
\<(\nabla_V A_{\tilde{\xi}})Y,Z\>=0.
$$

\noindent Using again that these equations are tensorial, we can
replace $Z$ by the vector field $X$. We obtain from
(\ref{eq:shape1}) and the Codazzi equation that
$$
\omega(\lambda\<V,e\>Y-\beta\<Y,X\>V)=0,
$$
which  by (\ref{eq:eqn1}) proves Fact \ref{claim} also in this
case.\vspace{1ex}

It follows from Theorem 5' in \cite{dt1} and the assumption that
$f$ is an embedding that $f$ is a composition $f=H\circ g$, where
$g\colon\,N^{p+n}\to\Q_c^{\,p+n+1}$ is an isometric immersion such
that $A^g_{\delta}=A_{\tilde{\xi}}$ for some unit normal vector
field $\delta$ of $g$, and $H\colon\, W\to\Q_c^{\,p+n+2}$ is an
\ii  of an open subset $W\subset \Q_c^{\,p+n+1}$ containing
$g(N^{p+n})$. Moreover, since $A^g_{\delta}=A_{\tilde{\xi}}$
satisfies  $A^g_{\delta}|_{{\mathcal V}_z}= \tilde{\gamma}_0\,id$
or $A^g_{\delta}|_{{\mathcal H}_z}= 0$ at any $z\in N^{p+n}$, it
follows that $g$ is of type $A$, and the conclusion follows from
Theorem \ref{thm:no}.\qed

\begin{remark}\label{re:typeb1} By Proposition \ref{pro:cij0},
if $f\colon\,N^{p+n}=L^p\times_{\rho} M^n\to\Q_c^{\,p+n+2}$, $n\ge
2$, is an isometric embedding   of type $B_1$ then it must satisfy
one of the conditions in Proposition~\ref{cor:typeb}-$(i)$  {\em
everywhere\/}. Moreover, for $n\geq 1$ we have that both
conditions hold simultaneously if and only if the \ii
$g\colon\,N^{p+n}\to\Q_c^{\,p+n+1}$ satisfies
$A^g_{\delta}|_{{\mathcal V}}= \tilde{\gamma}_0\,id$ and
$A^g_{\delta}|_{{\mathcal H}}= 0$. By Corollary \ref{le:triv},
this is the case if and only if \mbox{$g=\Psi\circ(h_1\times
h_2)$} with $h_1$ totally geodesic and  $h_2$ a local isometry,
where $\Psi\colon\,V^{p+1}\times_{\sigma} \Q_{\tilde{c}}^{\,n} \to
\Q_c^{\,p+n+1}$ is a warped product representation determined by
$(\Q_{\tilde{c}}^{\,n},\bar{z})$. Furthermore, if in addition
$\tilde{\gamma}_0=0$, that is, $g$ is totally geodesic, and
$\grad\,\rho$ has no zeros, then Corollary~\ref{cor:cyl} implies
that $h_1$ must be cylindrical with respect to $a$, where $-a$  is
the mean curvature vector of $\Q_{\tilde{c}}^{\,n}$  at $\bar{z}$
in either $\R^{p+n+1}$ or $\Oes^{\,p+n+2}$, according as $c=0$ or
$c\neq 0$. Conversely, if $h_1$ is totally geodesic and
cylindrical with respect to $a$ and $h_2$ is a local isometry,
then \mbox{$g=\Psi\circ(h_1\times h_2)$} is totally geodesic.
\end{remark}

  By using the main lemma in \cite{mo} or  Corollary \ref{thm:molzan}, according as $c=0$ or $c\neq 0$,  instead of Theorem \ref{thm:no}, we
obtain the following for isometric immersions of Riemannian
products.

\begin{corollary}\label{cor:rpb} Let
$f\colon\,L^p\times M^n\to \Q_c^{\,p+n+2}$, $p+n\ge 3$, be an
isometric embedding of type $B_1$ of a Riemannian product.
\begin{itemize}
\item[$(i)$] If $c=0$, then there exist an orthogonal decomposition $\R^{\,p+n+2}=\R^{\,p+k_1}\times\R^{\,n+k_2}$
with $k_1+k_2=1$ and \iis $h_1\colon\,L^p_0\to \R^{\,p+k_1}$,
\mbox{$h_2\colon M^n_0\to \R^{\,n+k_2}$} and $H\colon\,W\to
\R^{\,p+n+2}$ of an open subset $W\supset (h_1\times h_2)(U)$ of
$\R^{\,p+n+1}$ such that $f|_U=H\circ (h_1\times h_2)$ $($see the
diagram in Case $c=0$ - $(ii)$ of Theorem \ref{cor:prod2}$\,)$.

\item[$(ii)$]  If $c\neq 0$ there exist an
isometric embedding $\Phi\colon\,\Q_{c_1}^{\,p}\times
\Q_{c_2}^{\,n}\to \Q_c^{\,p+n+1}$ as an extrinsic Riemannian
product, local isometries $i_1\colon\, L^p\to \Q_{c_1}^{\,p}$ and
$i_2\colon\, M^n\to \Q_{c_2}^{\,n}$,  and an \ii $H\colon\,W\to
\Q_c^{\,p+n+2}$ of an open subset $W\supset \Phi\circ (i_1\times
i_2)(L^p\times M^n)$ of  $\Q_c^{\,p+n+1}$ such that $f=H\circ
\Phi\circ (i_1\times i_2)$ $($see the diagram in Case $c\neq 0$ -
$(ii)$ of Theorem~\ref{cor:prod2}$\,)$.
\end{itemize}
\end{corollary}

  We now consider  isometric
immersions $f\colon\,L^p\times_\rho M^n\to\Q_c^{\,p+n+2}$ of
type~$B_2$. In the following statement, in order not to have to
consider separately the cases $p=1$ and $p\geq 2$, we agree that
in the first case all information related to the splitting
$L^p=L^{p-1}\times_{\rho_1}I$ should be disregarded.

\begin{proposition}\label{le:B_5} Let
$f\colon\,L^p\times_\rho M^n\to\Q_c^{\,p+n+2}$ be an \ii of type
$B_2$ and assume that $n\geq 3$. Then locally we have:  $L^p$ and
$M^n$ split as warped products $L^p=L^{p-1}\times_{\rho_1}I$ and
$M^n =J\times_{\rho_2}M^{n-1}$ where $I,J\subset\R$ are open
intervals, and

$$
N^{p+n}=L^{p-1}\times_{\rho_1}
((I\times_{\rho_3}J)\times_{\bar{\rho}}M^{n-1}),
$$
where $\rho_1\in C^\infty (L^{p-1})$, $\rho_2\in C^\infty (J)$,
$\rho_3\in C^\infty (I)$ and $\bar{\rho}\in C^\infty (I\times J)$
satisfy
$$
\rho=(\rho_1\circ \pi_{L^{p-1}})(\rho_3\circ \pi_{I})
\;\;\mbox{and}\;\;\bar{\rho} =(\rho_3\circ \pi_{I})(\rho_2\circ
\pi_{J}),
$$
and there exist warped product representations
$$
\Psi_1\colon\,V^{p-1}\times_{\sigma_1} \Q_{\tilde{c}}^{\,n+3}\to
\Q_c^{\,p+n+2}\; \an\;
\Psi_2\colon\,W^{4}\times_{\sigma_2}\Q_{\bar{c}}^{\,n-1}
\to\Q_{\tilde{c}}^{\,n+3},
$$
isometries $i_1\colon\,L^{p-1} \to W^{p-1}\subset
V^{p-1}\subset\Q_c^{\,p-1}$ and  $i_2\colon\,M^{n-1} \to
W^{\,n-1}\subset \Q_{\bar{c}}^{\,n-1}$ onto open subsets, and an
\ii $g\colon\, I\times_{\rho_3}J\to W^{4}$ of type $B_2$ such that
$\bar{\rho} =\sigma_2\circ g$, $\rho_1=\sigma_1\circ i_1$
and~$f=\Psi_1\circ (i_1\times (\Psi_2\circ (g\times i_2)))$.
Moreover, $L^p$ has constant sectional curvature  $c$ if $p\geq 2$
$($see the diagram in Theorem \ref{thm:main}-$(iii))$.
\end{proposition}

\proof We have by (\ref{eq:shape11}), (\ref{eq:B_1}) and
(\ref{eq:B_2}) that \be\label{eq:sf2}
A_{\tilde{\xi}}Y=\tilde{\beta}_0\<Y,X\>X
\;\;\;\;\mbox{and}\;\;\;\; A_{\tilde{\xi}}V
=\tilde{b}V+(\tilde{\gamma}_0-\tilde{b})\<V,e\>e, \ee where
$\tilde{b}=\tilde{\beta}_0^{-1}\,\tilde{\delta}_0$. On the other
hand, we have that $A_\xi$ is given by (\ref{eq:shape1}). Thus,
for the relative nullity subspace $\Delta(z)$ at $z\in N^{p+n}$
there are two possibilities:
$$
\Delta(z)=\left\{\begin{array}{l} \{X\}^\perp\subset {\mathcal
H}\hspace*{17,2ex}\mbox{if}\;\;\tilde{\delta}_0(z)\neq 0,
\vspace*{2ex}\\
(\{X\}^\perp\subset {\mathcal H})\oplus (\{e\}^\perp \subset
{\mathcal V})\;\;\mbox{if}\;\;\tilde{\delta}_0(z)=0.
\end{array}\right.
$$

In the remaining of this proof the letters $T$ and $S$ will always
denote vector fields in $\Gamma(\{X\}^\perp\subset {\mathcal H} )$
and $\Gamma(\{e\}^\perp \subset {\mathcal V})$, respectively. We
also denote by $(A_\delta,u,v,w)$ taking the $w$-component of the
Codazzi equation for $A_\delta$ and the vectors $u,v$.

    We first prove that if $p\geq 2$ then $L^p$ splits locally as $L^p=L^{\,p-1}\times_{\rho_1} I$,
where $I$ is an open interval, and that $\rho=(\rho_1\circ
\pi_{L^{\,p-1}})(\rho_3\circ \pi_{I})$ for some functions
$\rho_1\in C^\infty (L^{\,p-1})$ and $\rho_3\in C^\infty (I)$. We
point out that this fact also holds if $n=1,2$. As a first step,
we show that the vector field $X$ is the lift of a vector field
$\tilde{X}\in \Gamma(TL)$. For that, we must prove that
\be\label{eq:xlift} \nabla_V X=-\<X,\eta\>V. \ee Notice that
$(A_{\tilde{\xi}},X,V,T)$ reads
$$
\<\nabla_X A_{\tilde{\xi}}V- A_{\tilde{\xi}}\nabla_XV-
\nabla_VA_{\tilde{\xi}}X+A_{\tilde{\xi}}\nabla_VX,T\>
=\<\omega(V)A_\xi X-\omega(X)A_\xi V,T\>.
$$
Then, by means of (\ref{eq:shape1}) and (\ref{eq:sf2}) we obtain
\be\label{eq:t5iv} \<\nabla_V X,T\>=0. \ee On the other hand,   by
means of (\ref{eq:c4}), \be\label{eq:t5x} \<\nabla_V
X,W\>=-\<\nabla_V W,X\>=-\<V,W\>\<X,\eta\>. \ee Using also that
$\<\nabla_V X,X\>=0$, for $X$ has unit length, we obtain
(\ref{eq:xlift}) from (\ref{eq:t5iv}) and (\ref{eq:t5x}). Notice
that $\tilde{X}$ has unit length, because
$$\<\tilde{X},\tilde{X}\>_L\circ\pi_L=\<{\pi_L}_* X, {\pi_L}_*
X\>_L=\<X,X\>_N=1.$$

    We show next that the distribution
$\{\tilde{X}\}^\perp$ is totally geodesic in $L^{\,p}$. In effect,
for any $\tilde{T}_1, \tilde{T}_2\in \Gamma(\{\tilde{X}\}^\perp)$
we have
$$
\<\nabla^L_{\tilde{T}_1} \tilde{T}_2,\tilde{X}\>_L\circ\pi_L
=\<{\pi_L}_*\nabla_ {T_1} T_2,{\pi_L}_*X\>_L =\<\nabla_{T_1}
T_2,X\>_N=0,
$$
where the last equality follows from the fact that the relative
nullity distribution $\Delta$ is totally geodesic and $T_1, T_2\in
\Gamma(\Delta)$, whereas $X\in \Gamma(\Delta^\perp)$.

   Our next step is to prove that the vector field $\zeta=\eta-\<X,\eta\>X\in \Gamma(\{X\}^\perp)$ is the
   lift of a vector field $\tilde{\zeta}\in \Gamma(\{\tilde{X}\}^\perp)$,  i.e., $\nabla_V \zeta=-\<\zeta,\eta\>V$.  This follows from
$$
\<\nabla_V \zeta,W\>= -\<\zeta,\nabla_V
W\>=-\<V,W\>\<\eta,\zeta\>,\;\;\;
  \<\nabla_V \zeta,X\>=-\<\zeta,\nabla_V X\>=0,
$$
where we have used (\ref{eq:xlift}), and from
$$
\<\nabla_V \zeta,T\>=\<\nabla_V \eta,T\>
-V(\<X,\eta\>)\<X,T\>-\<X,\eta\>\<\nabla_V X,T\>=0,
$$
where we have used (\ref{eq:xlift}) for the last term and that
${\mathcal V}$ is spherical for the first term.

     Our final step is to show that the distribution $\{\tilde{X}\}$ is spherical with mean curvature vector $\tilde{\zeta}$. We have from $(A_{\xi},X,e,T)$  and $\lambda\neq 0$ that
$\<\nabla_X X,T\>=\<\nabla_e e,T\>=\<\eta,T\>$, and hence
$$
\nabla^L_{\tilde{X}} \tilde{X}\circ \pi_L={\pi_L}_*\nabla_X
X={\pi_L}_*\zeta=\tilde{\zeta}.
$$
On the other hand, we obtain from (\ref{eq:ge11}) for $Y=T$ and
$V=W\neq 0$ that
$$
\<\nabla_X \eta, T\>=\<X,\eta\>\<T,\eta\>.
$$
Thus,
$$
\<\nabla_X \zeta, T\>=\<\nabla_X \eta, T\> -\<X,\eta\>\<\nabla_X
X,T\>=0,
$$
and therefore,
$$
\<\nabla^L_{\tilde{X}} \tilde{\zeta}, \tilde{T}\>_L \circ \pi_L
=\<{\pi_L}_*\nabla_X \zeta, {\pi_L}_* T\>_L =\<\nabla_X \zeta,
T\>_N=0 \;\;\mbox{for }\;\;\;\tilde{T}\in
\Gamma(\{\tilde{X}\}^\perp),
$$
which completes the proof of the step.

    By Theorem \ref{thm:hi}, we have that locally $L^p$ splits as $L^p=L^{\,p-1}\times_{\rho_1} I$,
where $I$ is an open interval and $\tilde{\zeta}=-\grad\, \log
(\rho_1\circ \pi_{L^{\,p-1}})$. In particular, the lift $\zeta$ of
$\tilde{\zeta}$ to $N^{p+n}$ is
$$
\zeta=-\grad\,\log (\rho_1\circ\pi_{L^{\,p-1}}\circ\pi_{L^{\,p}}).
$$
Since we also have $\eta=-\grad\log (\rho\circ \pi_L)$, we obtain
$\<X,\eta\>X=\eta-\zeta=-\grad\log (\hat{\rho}\circ \pi_L)$ with
$\hat{\rho}=\rho(\rho_1\circ\pi_{L^{p-1}})^{-1} \in
C^{\infty}(L^p)$. Moreover, since
$$
(\tilde{T}(\log \hat{\rho}))\circ \pi_{L^p}=T(\log(\hat{\rho}\circ
\pi_{L^p}))=-\<X,\eta\>\<X,T\>=0,
$$
it follows that there exists $\rho_3\in C^\infty(I)$ such that
$\hat{\rho}=\rho_3\circ \pi_I$. \vspace{1ex}

    Let us now prove that locally $M^n$ also splits as
$M^n=J\times_{\rho_2}M^{n-1}$, where $J$ is an open interval.
First, we obtain from $(A_{\tilde{\xi}},Y,e,S)$ and $\tilde b\neq
\tilde{\gamma}_0$ that  $\<\nabla_Y e,S\>=0$. Since also
$\<\nabla_Y e,e\>=0$, for $e$ has unit length, and $\<\nabla_Y
e,Z\>=-\<\nabla_Y Z,e\>=0$,  we have \be\label{eq:nye} \nabla_Y
e=0. \ee It follows that
$$\nabla_Y(\rho\circ \pi_L)e=Y(\rho\circ \pi_L)e=-\<Y,\eta\>(\rho\circ \pi_L)e.$$
This implies that $(\rho\circ \pi_L) e$ is the lift of a vector
field $\tilde{e}\in\Gamma(M^n)$. Notice that $\tilde{e}$ is a unit
vector field, for
$$
\<\tilde{e},\tilde{e}\>_M\circ \pi_M=\<{\pi_M}_*(\rho e),
{\pi_M}_*(\rho e)\>_M=\rho^{-2}\<\rho e,\rho e\>_N=1.
$$
Thus, in order to show that locally $M^n$ splits as claimed,  by
Theorem \ref{thm:hi} it suffices to prove that the distribution
$\{\tilde{e}\}$ is totally geodesic and that $\{\tilde{e}\}^\perp$
is spherical.

 We obtain from $(A_\xi,X,e,S)$ and (\ref{eq:nye}) that
\be\label{eq:t5i} \<\nabla_e e,S\>=0. \ee In particular, it
follows that $\nabla_e e=\eta$. We conclude that the distribution
$\{\tilde{e}\}$ is totally geodesic from
$$
\<\nabla^M_{\tilde{e}} \tilde{e},\tilde{S}\>_M \circ\pi_M
=\<{\pi_M}_*\nabla_{\rho e} \rho e,{\pi_M}_*S\>_M =\<\nabla_e
e,S\>_N=0
$$
for any $\tilde{S}\in \Gamma(\{\tilde{e}\}^\perp)$. On the other
hand, we obtain from $(A_{\tilde{\xi}},S_1,e,S_2)$ that
$$
\<\nabla_{S_1} S_2,e\>=\va\<S_1,S_2\>,
$$
where $\va=(\tilde{b}-\tilde{\gamma}_0)^{-1}e(\tilde b)$. Thus,
the distribution $\{e\}^\perp\subset \V$ is totally umbilical.
Moreover,
$$
\begin{array}{l}
\<\nabla^M_{\tilde{S}_1} \tilde{S}_2,\tilde{e}\>_M\circ\pi_M
=\<{\pi_M}_*\nabla_{S_1} S_2, {\pi_M}_*\rho e\>_M
=\rho^{-1}\<\nabla_{S_1} S_2, e\>_N\\\hspace{18.7ex}
=\rho^{-1}\va\<S_1,S_2\>_N=\rho\va\<\tilde{S}_1,\tilde{S}_2\>_M\circ\pi_M.
\end{array}
$$
The preceding equality implies that there exists $\tilde{\va}\in
C^\infty(M^n)$ such that \mbox{$\rho\va=\tilde{\va}\circ \pi_{M}$}
and that the distribution $\{\tilde{e}\}^\perp$ is totally
umbilical with mean curvature normal $\tilde{\va}\tilde{e}$.  {\em
In particular, if $\tilde{b}$, or equivalently,
$\tilde{\delta}_0$, vanishes on an open subset $L_0^p\times
M_0^n\subset N^{p+n}$, then $\{\tilde{e}\}^\perp$ is a totally
geodesic distribution in $M_0^n$\/}. In the general case, in order
to show that $\{\tilde{e}\}^\perp$ is spherical, it remains to
prove that $\tilde{S}(\tilde{\va})=0$ or, equivalently,
that~$S(\va)=~0$.   First, using that ${\mathcal V}$ is umbilical
and invariant by $A_{\tilde{\xi}}$, and that $\lambda\neq 0$, we
obtain from $(A_{\tilde{\xi}},e,S,X)$ that
$\nabla^\perp_S\tilde{\xi}=0$. Now, choosing linearly independent
sections $S_1,S_2\in \Gamma(\{e\}^\perp)$ we obtain from
$(A_{\tilde \xi},S_1,S_2,S_1)$ that $S(\tilde{b})=0$.  {\em  We
point out that the assumption that $n\geq 3$ is only used here. In
particular, if  $\tilde{\delta}_0$ is everywhere vanishing then it
is enough to assume that $n\geq 2$.\/} Using (\ref{eq:t5i}) we
obtain from $(A_{\tilde{\xi}},e,S,e)$ that
$S(\tilde{\gamma}_0)=0$. Since $\nabla_e S\in \{e\}^\perp$ and
$\nabla_S e\in \{e\}^\perp$, as follows from (\ref{eq:t5i}) and
the fact that ${\mathcal V}$ is totally umbilical, then
$$
Se(\tilde{b})=eS(\tilde{b})+\nabla_e S(\tilde{b}) -\nabla_S
e(\tilde{b})=0,
$$
and hence $S(\va)=0$.  Therefore, locally $N^{p+n}$  splits as
$$
N^{p+n}=L^{p-1}\times_{\rho_1}M^{n+1}\;\;\mbox{with}\;\;
M^{n+1}:=M^2\times_{\bar{\rho}}M^{n-1}\;\mbox{and}\;
M^2:=I\times_{\rho_3}J,
$$
where $J\subset \R$ is an open interval and $\bar{\rho}=
(\rho_3\circ \pi_{I})(\rho_2\circ \pi_{J})$. Notice that $f$ is of
type $A$ with respect to this decomposition of $N^{p+n}$. We claim
that there exist a warped product representation
$\Psi_1\colon\,V^{p-1} \times_{\sigma_1}\Q_{\tilde{c}}^{\,n+3} \to
\Q_c^{\,p+n+2}$, an isometric immersion
$\tilde{G}\colon\,M^{n+1}\to \Q_{\tilde{c}}^{\,n+3}$ and a local
isometry $i_1\colon\,L^{p-1}\to V^{p-1}$  such that
$\rho_1=\sigma_1\circ i_1$ and $f=\Psi_1\circ (i_1\times
\tilde{G})$.

 Fix $\bar{y}\in L^{p-1}$ with $\rho_1(\bar{y})=1$ and let
 $i_{\bar{y}}\colon\,M^{n+1}\to N^{n+p}$ be the (isometric) inclusion of $M^{n+1}$ into
 $N^{n+p}$ as a leaf of the vertical subbundle $\bar{{\mathcal V}}$ according to the latter
 decomposition of $N^{n+p}$. Define \mbox{$G\colon\,M^{n+1}\to \Q_c^{\,p+n+2}$} by
 $G=f\circ i_{\bar{y}}$.  By Corollary
\ref{cor:shull2}-$(ii)$, in order to show that the spherical hull
of $G$ has dimension $n+3$, it suffices to prove that for no point
$z\in N^{n+p}$ there exists a unit vector $\bar{\xi}\in T_z^\perp
N$ such that $A_{\bar{\xi}}|_{\bar{{\mathcal
V}}_z}\colon\,{\bar{{\mathcal V}}_z}\to {\bar{{\mathcal V}}_z}$ is
a multiple of the identity tensor.  Write $\bar{\xi}=\cos
\theta\xi+\sin \theta\tilde{\xi}$. Then $\<A_{\bar{\xi}}X,e\>=0$
and $\<A_{\bar{\xi}}X,X\>=\<A_{\bar{\xi}}e,e\>
=\<A_{\bar{\xi}}S,S\>$ for any unit vector  $S\in
(\{e\}_z^\perp\subset {\mathcal V}_z)$ if and only if
$$
\lambda\cos \theta=0\;\;\;\mbox{and}\;\;\;\tilde{\beta}_0\sin
\theta+\beta\cos
\theta=\gamma\cos\theta+\tilde{\gamma}_0\sin\theta=\tilde{b} \sin
\theta.
$$
Since $\lambda\neq 0$, we obtain that
$\tilde{\gamma}_0=\tilde{b}$,
 a contradiction to the fact that $\beta\gamma-\lambda^2\neq 0$.
 Our claim then follows from Theorem
\ref{thm:no} by letting $\Q_{\tilde{c}}^{\,n+3}$ be the spherical
hull of $G$ and defining $\tilde{G}\colon\,M^{n+1}\to
\Q_{\tilde{c}}^{\,n+3}$  by $G=j\circ \tilde{G}$, where $j$ is the
inclusion of $\Q_{\tilde{c}}^{\,n+3}$ into $\Q_c^{\,p+n+2}$.

    We now study the isometric immersion $\tilde{G}\colon\,M^{n+1}\to
\Q_{\tilde{c}}^{\,n+3}$. First observe that the second fundamental
form of $G=j\circ \tilde{G}$ is given by
$$\alpha_G(V,W)=\alpha_f({i_{\bar{y}}}_*V,{i_{\bar{y}}}_*W)+\<V,W\>f_*(\bar{\eta}\circ
i_{\bar{y}})\;\;\;\mbox{for all }\;\;V,W\in \Gamma(TM^{n+1}),$$
where $\bar{\eta}$ denotes the mean curvature normal of
$\bar{\V}$. Let $\hat{\V}$ denote the vertical subbundle of
$TM^{n+1}$ according to the decomposition
$M^{n+1}=M^2\times_{\bar{\rho}}M^{n-1}$. Using that
$A^f_\xi|_{{i_{\bar{y}}}_*\hat{\V}}=0$ and that
$A^f_{\tilde{\xi}}|_{{i_{\bar{y}}}_*\hat{\V}}=\tilde{b}\,id$,
where $id$ denotes the identity tensor, it follows that
$$\alpha_G(V,\bar{V})=\<V,\bar{V}\>((\tilde{b}\circ
i_{\bar{y}})(\tilde{\xi}\circ i_{\bar{y}})+f_*(\bar{\eta}\circ
i_{\bar{y}}))\;\;\;\mbox{for all }\;\;V\in
\Gamma(TM^{n+1}),\;\;\bar{V}\in \hat{\V}.$$
 Therefore $\hat{\eta}=(\tilde{b}\circ
i_{\bar{y}})(\tilde{\xi}\circ i_{\bar{y}})+f_*(\bar{\eta}\circ
i_{\bar{y}})$ is a principal curvature normal of $G$ and
$\hat{\V}$ is contained in the corresponding eigendistribution.
Since $j$ is umbilical, it follows that
$\hat{\eta}_{T\Q_{\tilde{c}}^{\,n+3}}$ is a principal curvature
normal of $\tilde{G}$ with the same eigendistribution as
$\hat{\eta}$.  By means of Corollary \ref{cor:shull2}-$(i)$ and
Theorem~\ref{thm:no}, we conclude that there exist a warped
product representation
$\Psi_2\colon\,W^{4}\times_{\sigma_2}\Q_{\bar{c}}^{\,n-1} \to
\Q_{\tilde{c}}^{\,n+3}$, an \ii \mbox{$g\colon\,M^2\to W^{4}$} and
a  local isometry $i_2\colon\,M^{n-1}\to \Q_{\bar{c}}^{\,n-1}$
such that $\bar{\rho}=\sigma_2\circ g$ and $\tilde{G}=\Psi_2\circ
(g \times i_2)$. Moreover, since the \sff of $g$ is determined by
the restriction of $\alpha_G$ to the horizontal subbundle
$\hat{\H}$ of $TM^{n+1}$ according to the decomposition
$M^{n+1}=M^2\times_{\bar{\rho}}M^{n-1}$, and hence by the
restriction of $\alpha_f$ to $\spa\{X,e\}$ (see formula
(\ref{eq:sff2})), it follows that $g$ is of type $B_2$.

Finally, since $\{X\}^\perp\subset\H$ is contained in $\Delta$, we
obtain that the curvature-like tensor $C$ defined in Proposition
\ref{prop:symm} satisfies
$$
C(Y_1,Y_2,Y_3,Y_4)=0\;\;\; \mbox{for
all}\;\;Y_1,Y_2,Y_3,Y_4\in\Gamma(\H).
$$
The last assertion then follows from the fact that, for a fixed
$\bar{x}\in M^n$, the inclusion $i_{\bar{x}}\colon\, L^p\to
N^{p+n}$ given by $i_{\bar{x}}(y)=(y,\bar{x})$ is a totally
geodesic isometric immersion.\qed

\begin{corollary}\label{cor:prodb2} Let $N^{p+n}=L^p\times M^n$ be a Riemannian product of
dimension $p+n\geq 3$. Then there exists no \ii
$f\colon\,N^{p+n}\to\Q_c^{\,p+n+2}$ of type $B_2$ if $c\neq 0$.
\end{corollary}

\proof We may assume $p\ge 2$. It follows from Proposition
\ref{le:B_5} that locally $L^p$ splits  as a Riemannian product
$L^p=L^{p-1}\times I$. The statement now follows from the fact
that $L^p$ has constant sectional curvature $c$.\qed\vspace{1ex}

  In order to complete the classification of \iis $f\colon\,L^p\times M^n\to\Q_c^{\,p+n+2}$ of type $B$ of Riemannian products of dimension $p+n\geq 3$, it remains
to determine those that are of type $B_2$ for $c=0$. Observe that
for such \iis equation (\ref{eq:B_2}) in Proposition
\ref{cor:typeb}
 holds with $\tilde{\delta}_0=0$ (see Corollary \ref{cor:typebprod}).
In the following result we solve the more general problem of
classifying \iis of type $B_2$ of {\it warped\/} products
satisfying that condition.

\begin{corollary}\label{le:b6} Let
$f\colon\,N^{p+n}=L^p\times_\rho M^n\to\Q_c^{\,p+n+2}$, $n\geq 2$,
 be an \ii of type $B_2$ for which $\tilde{\delta}_0$ vanishes everywhere. Then $c\leq 0$,
 $N^{p+n}$ has constant sectional curvature $c$ and one of the
following holds locally:
\begin{itemize}
\item[$(i)$] If $c=0$ then $L^p$ and $M^n$
split as Riemannian products $L^p=L^{p-1}\times I$ and
$M^n=J\times M^{n-1}$, where $I,J\subset \R$ are open intervals,
and there exist isometries $i_1\colon\,L^{p-1}\to U\subset
\R^{p-1}$ and $i_2\colon\,M^{n-1}\to V\subset \R^{n-1}$ onto open
subsets and an \ii $g\colon\,I\times J\to \R^4$ such that
$f=i_1\times g\times i_2$.

\bigskip
\vspace{1ex}
\begin{picture}(150,84)\hspace{-13ex}
\put(100,10){$L^{p-1}\times\,(I\times J) \times M^{n-1} =N^{p+n}$}
\put(123,28){$i_1$} \put(232,45){$f=i_1\times g\times i_2$}
\put(158,45){$g$} \put(194,28){$i_2$}
\put(118,22){\vector(0,1){20}} \put(226,22){\vector(0,1){50}}
\put(153,22){\vector(0,1){50}} \put(188,22){\vector(0,1){20}}
\put(115,65){$\cup$} \put(115,50){$U$} \put(185,65){$\cup$}
\put(185,50){$V$} \put(115,80){$\R^{p-1}\times\, \R^4\;
\times\;\,\R^{n-1}=\,\R^{p+n+2}$}
\end{picture}

\item[$(ii)$] If $c<0$ then $L^p$ splits as a warped product
$L^p=L^{p-1}\times_{\rho_1}I$, $M^n$ splits as a Riemannian
product $M^n=J\times M^{n-1}$, where $I,J\subset \R$ are open
intervals and $\rho=\rho_1\circ \pi_{L^{p-1}}$, and there exist a
warped product  representation $\Psi\colon\,
{V^{p-1}}\times_\sigma\R^{n+3} \to \Q_c^{\,p+n+2}$, isometries
$i_1\colon\,L^{p-1}\to U\subset V^{\,p-1}$ and
\mbox{$i_2\colon\,M^{n-1}\to W\subset\R^{n-1}$} onto open subsets,
and an \ii \mbox{$g\colon\,I\times J\to \R^4$} such that
$f=\Psi\circ (i_1\times (g\times i_2))$ and $\rho_1=\sigma\circ
i_1$.
\end{itemize}
\bigskip
\begin{picture}(150,84)\hspace{-12ex}
\put(82,10){$L^{p-1}\times_{\rho_1} ((I\times J)\times M^{n-1})$}
\put(148,40){$g$} \put(340,10){$\Q_c^{\,p+n+2}$}
\put(222,-4){$f=\Psi\circ (i_1\times (g\times i_2))$}
\put(102,28){$i_1$} \put(194,28){$i_2$}
\put(143,22){\vector(0,1){40}} \put(96,22){\vector(0,1){20}}
\put(178,22){\vector(0,1){20}} \put(93,58){$\cup$}
\put(93,46){$U$} \put(176,58){$\cup$} \put(174,46){$W$}
\put(92,70){$V^{p-1}\times_\sigma \;(\R^4\,\times\;\R^{n-1})$}
\put(255,26){$\circlearrowright$} \put(262,52){$\Psi$}
\put(203,65){\vector(3,-1){134}} \put(212,12){\vector(1,0){120}}
\end{picture}
\bigskip
\end{corollary}

\proof First observe  that  for the statement of Proposition
\ref{le:B_5} to hold in this case it is enough to require that
$n\ge 2$, as follows from the second italicized text in its proof.
In order to prove that $N^{p+n}$ has constant  sectional curvature
$c$, we must show that the curvature-like tensor $C$ defined in
Proposition~\ref{prop:symm} vanishes identically. Since the
relative nullity distribution of $f$ is
$\Delta=(\spa\{X,e\})^\perp$, we have that $C(E_1,E_2,E_3,E_4)=~0$
whenever two of the vectors $E_1,E_2,E_3,E_4$ belong to
$(\spa\{X,e\})^\perp$. Thus it remains to show that
$C(X,e,e,X)=0$, because $C$ is a curvature-like tensor. But this
follows from (\ref{eq:shape1}), (\ref{eq:sf2}) and the assumption
that $\tilde{\delta}_0=0$.

      We have from Proposition \ref{le:B_5} that locally $L^p$ and $M^n$ split  as
warped products $L^p=L^{p-1}\times_{\rho_1}I$, $M^n
=J\times_{\rho_2}M^{n-1}$,  and $N^{p+n}=L^{p-1}\times_{\rho_1}
((I\times_{\rho_3}J)\times_{\bar{\rho}}M^{n-1}), $ where
$I,J\subset\R$ are open intervals and $\rho_1\in C^\infty
(L^{p-1})$, $\rho_2\in C^\infty (J)$, $\rho_3\in C^\infty (I)$ and
$\bar{\rho}\in C^\infty (I\times J)$ satisfy
$$
\rho=(\rho_1\circ \pi_{L^{p-1}})(\rho_3\circ \pi_{I})
\;\;\mbox{and}\;\;\bar{\rho} =(\rho_3\circ \pi_{I})(\rho_2\circ
\pi_{J}).
$$
But now $ \<X,\eta\>=\<\nabla_S S,X\>=0, $ because $\Delta$ is
totally geodesic.  Hence we may assume that $\rho_3=1$ (recall the
proof of Proposition \ref{le:B_5}), and therefore
$\rho=\rho_1\circ \pi_{L^{p-1}}$. On the other hand, by the first
italicized text in the proof of Proposition \ref{le:B_5}, the
distribution $\{\tilde{e}\}^\perp$ in $M^n$ is now totally
geodesic, and hence $\rho_2=1$, which implies that also
$\bar{\rho}=1$. Summing things up, we have
$$
N^{p+n}=L^{p-1}\times_{\rho_1}M^{n+1},\;\;
\mbox{with}\;\;\;M^{n+1}:=M^2\times M^{n-1}\;\mbox{and}\;\;
M^2:=I\times J.
$$

    Now, for a fixed point $\bar{y}\in L^p$ with $\rho(\bar{y})=1$,
    let $i_{\bar{y}}\colon\,M^{n}\to N^{p+n}$  denote the (isometric) inclusion of
    $M^{n}$ into $N^{p+n}$ as a leaf of ${\V}$.  The second fundamental form of
    $i_{\bar{y}}$ is $\alpha_{i_{\bar{y}}}(V,W)=\<V,W\>(\eta\circ
    i_{\bar{y}})$ for all $V, W\in \Gamma(TM^n)$, where
    $\eta=-\grad \log(\rho\circ \pi_L)$ is the mean curvature
    normal of $\V$. Since $N^{p+n}$ has constant sectional curvature $c$, it follows from
    the Gauss equation for
    $i_{\bar{y}}$  that $M^n$ has constant sectional curvature
    $c+\|\eta\circ i_{\bar{y}}\|^2=c+\|\grad \log
    \rho(\bar{y})\|^2$. We conclude from the fact that $M^n=J\times M^{n-1}$ is a Riemannian
    product that it must be flat, hence $M^{n-1}$ must be flat when $n\geq 3$ and
    $c+\|\grad \log \rho(\bar{y})\|^2=0$. Now choose any other point $y^*\in L^p$ and modify
    the warped product representation of $N^{n+p}$ so that the modified warping function
    $\rho^*$ satisfies $\rho^*(y^*)=1$. By this modification
    $\eta=-\grad \log(\rho\circ \pi_L)$ does not change. Therefore the preceding argument yields
    $\|\grad \log \rho({y}^*)\|^2=-c$. \vspace{1ex}\\We now distinguish the
two possible cases:\vspace{1ex}\\
{\em Case $c=0$.\/} Here $\eta=-\grad (\log \rho\circ \pi_{L^p})$
vanishes, hence $\rho=1$, consequently also $\rho_1=1$ and
therefore $N^{p+n}=L^{p-1}\times M^2 \times M^{n-1}$, with
$L^{p-1}$ and $M^{n-1}$ flat and $M^2=I\times J$.  Using that
$\Delta=(\spa\{X,e\})^\perp$ is the relative nullity distribution
of $f$, the main lemma in \cite{mo} implies that $f$ splits as
$$
f=i_1\times g\times  i_2\colon\,L^{p-1}\times M^2\times M^{n-1}
\to \R^{p-1}\times\R^4\times\R^{n-1}=\R^{p+n+2},
$$
where $i_1\colon\,L^{p-1}\to U\subset \R^{p-1}$ and
$i_2\colon\,M^{n-1}\to V\subset \R^{n-1}$ are isometries onto open
subsets, and $g\colon\,M^2\to \R^4$ is an
isometric immersion.\vspace{1ex}\\
{\em Case $c<0$.\/} For $G=f\circ i_{\bar{y}}$ as in the proof of
Proposition \ref{le:B_5}, we have from Corollary
\ref{cor:shull2}-$(ii)$ that its  spherical hull
$\Q_{\tilde{c}}^{n+3}$ has constant sectional curvature
$\tilde{c}=c+\|\grad \log \rho(\bar{y})\|^2=0$. Therefore
$\Q_{\tilde{c}}^{n+3}$ is a horosphere
$\R^{n+3}\subset\Q_c^{\,p+n+2}$. Let $\tilde{G}\colon\,M^{n+1}\to
\R^{n+3}=\Q_{\tilde{c}}^{\,n+3}$ be such that $G=j\circ
\tilde{G}$, where $j$ denotes the inclusion of
$\Q_{\tilde{c}}^{\,n+3}$ into $\Q_c^{\,p+n+2}$. Using that the
vertical subbundle of $M^{n+1}$ correspondent to the splitting
$M^{n+1}=M^2\times M^{n-1}$ is contained in the relative nullity
distribution of $\tilde{G}$, the conclusion now follows from the
main lemma in \cite{mo} applied to $\tilde{G}\colon\,M^2\times
M^{n-1}\to \R^{n+3}$.\qed

\begin{corollary}\label{escolio} Let $N^{p+n}=L^p\times M^n$ be a Riemannian product of dimension $p+n\geq 3$.
Then any \ii $f\colon\,N^{p+n}\to\R^{p+n+2}$ of type $B_2$ is
locally given as in Corollary \ref{le:b6}-$(i)$.

\end{corollary}

\begin{remark}\label{bdt} By making use of global
arguments from \cite{am}, a complete description of the possible
cases in which an \ii $f\colon\,L^p\times M^n\to \R^{\,p+n+2}$,
$p\ge 2$ and $n\ge 2$, of a Riemannian product of complete nonflat
Riemannian manifolds may fail locally to be a product of isometric
immersions was given in \cite{bdt}. Namely, it was shown therein
that there exists an open dense subset of $L^p\times M^n$ each of
whose points lies in an open product neighborhood $U_0=L_0^p\times
M_0^n$ restricted to which $f$ is either $(i)$ a product of
isometric immersions, $(ii)$ an isometric immersion of type $B_2$
given as in  Corollary \ref{le:b6} - $(i)$, or $(iii)$ an
isometric immersion of type $B_1$ of the following special type:
either $L_0^p$ or $M_0^n$, say, the latter, splits as
$M_0^n=I\times \R^{n-1}$, the manifold $L_0^p$ is free of flat
points and $f|_{U_0}$ splits as
$$
f|_{U_0}=F\times id\colon\,(L_0^p\times I) \times \R^{n-1}\to
\R^{p+3}\times \R^{n-1}=\R^{p+n+2}.
$$
Moreover,  $F\colon\,L_0^p\times I\to \R^{p+3}$ is a composition
$F=H\circ \tilde{F}$, where
$$
\tilde{F}=G\times i\colon\,L_0^p\times I\to \R^{p+1}\times
\R=\R^{p+2}
$$
is a cylinder over a hypersurface  $G\colon\,L_0^p\to\R^{p+1}$,
and $H\colon\,W\to \R^{p+3}$ is an \ii of an open subset $W\supset
\tilde{F}(L_0^p\times I)$ of $\R^{p+2}$.

We take the opportunity to point out that the main theorem in
\cite{bdt} misses the informations that $L_0^p$ is free of flat
points and that $\tilde{F}$ is a cylinder $\tilde{F}=G\times i$,
which follow  from Corollary \ref{cor:prod21}.
\end{remark}

\section{Immersions of type C}

  The aim of this section is to prove the following
pointwise result for \iis  of type $C$.

\begin{proposition}\label{casec} Let
$f\colon\,L^p\times_\rho M^n\to\Q_c^{\,p+n+2}$ be an isometric
immersion of a warped product. If $f$ is of type $C$ at $z\in
N^{p+n}=L^p\times_\rho M^n$ and $n\geq 3$, then $N^{p+n}$ has
constant sectional curvature $c$ at $z$.
\end{proposition}

\proof We must prove that the curvature like tensor $C$ on $T_zN$
defined in Proposition~\ref{prop:symm} vanishes identically. Two
possible cases may occur: \vspace{1ex}

\noindent Case $1$.  For every $X\in {\mathcal H}_z$ the linear
map $B_X$ defined in (\ref{eq:bx}) satisfies $\rank B_X\leq
1$.\vspace{1ex}

\noindent Case $2$.  There exists $X\in {\mathcal H}_z$ such that
$\rank B_X=2$.\vspace{1ex}

We first prove the following facts. \vspace{1ex}

\noindent $(i)$  In Case $1$, for any $X\in\H_z$ with $\rank
B_X=1$ we have \be \label{eq:aev} \a(E,V)=0\;\fall V\in\D(X)=\ker
B_X\an  E\in T_zN, \ee that is, $\D(X)$ is contained in the
relative nullity subspace $\Delta$ of $f$ at $z$.\vspace{1ex}

\noindent $(ii)$ In Case $2$, condition (\ref{eq:aev}) is true for
any $X\in\H_z$ such that $\rank B_X=2$. \vspace{.5ex}

\noindent {\em For $(i)$.\/} Here  Lemma \ref{prop:rankone}
applies and we have by (\ref{eq:r1a}) that $\D(Y)=\{e\}^\perp$ for
any $Y\in\H_z$ with $\rank B_Y=1$. In particular, this shows that
(\ref{eq:aev}) is satisfied for any $E\in {\mathcal H}_z$. On the
other hand, (\ref{eq:r1b}) yields $\a(V,W)=\<V,e\>\<W,e\>\a(e,e)$,
and hence (\ref{eq:aev}) also holds for any
$E\in {\mathcal V}_z$.\vspace{1ex}\\
\noindent {\em For $(ii)$.\/} If $\rank B_X=2$, i.e.,
$B_X({\mathcal V}_z)=T_z^\perp N$, then (\ref{eq:aev}) is
equivalent to
$$
\<\a(E,V),\a(X,W)\>=0
$$
for all $V\in \D(X)$, $W\in\V_z$ and $E\in T_zN$. But this follows
from (\ref{eq:ge33}) for $E\in\H_z$ and from (\ref{eq:ge44}) for
$E\in\V_z$. \vspace{1ex}

We now prove that \be\label{c} C(E_1,V_1,V_2,E_2)=0\;\;\;\mbox{for
all}\;\;E_1,E_2\in T_zN \an V_1,V_2\in\V_z. \ee We obtain from
(\ref{eq:ge44}) that (\ref{c}) holds whenever one of the vectors
$E_1,E_2$ lies in $\H_z$ and the other in $\V_z$. On the other
hand, if we are in Case $1$ (resp., Case~$2$) and $X\in\H_z$
satisfies $\rank B_X=1$ (resp., $\rank B_X=2$), then it follows
from (\ref{eq:aev}) that (\ref{c}) is also satisfied if any of the
vectors $E_1,E_2,V_1$ or $V_2$ belongs to $\D(X)$. In particular,
$C(E_1,V,V ,E_2)=0$ holds for any $V\in\D(X)$. Notice that
$\D(X)\neq\{0\}$ by our assumption that $n\ge 3$. Applying
(\ref{eq:ge11}) for $0\neq W=V\in \D(X)$ yields
\be\label{eq:point}\nabla_Y\eta-\<Y,\eta\>\eta-cY=0\;\; \mbox{for
any}\;\;Y\in\H_z. \ee

 {\em Notice that if $N^{p+n}$ is a Riemannian product
then this implies that $c=0$. Moreover, since for $c=0$ this
equation holds automatically for Riemannian products, in this case
it is enough to assume that either $n\ge 2$ or $p\geq 2$.\/
}\vspace{.5ex}

Therefore, (\ref{c}) is satisfied for all $E_1,E_2\in\H_z$ and
$V_1,V_2\in\V_z$. This completes the proof of (\ref{c}) in Case
$1$ and shows that in Case $2$ it remains to prove that (\ref{c})
is satisfied if $E_1,E_2,V_1,V_2$ all belong to the
two-dimensional subspace $\D(X)^\perp$. Since $C$ is a curvature
like tensor, this will follow once we prove the existence of an
orthonormal basis $\{e_1,e_2\}$ of $\D(X)^\perp$ such that
\be\label{flat} C(e_1,e_2,e_2,e_1)=0. \ee

  In order to prove (\ref{flat}) take an  orthonormal basis
$e_1,e_2$ of $\D(X)^\perp$ such that $e_1$ is one of the two
points on the unit circle $S^1$ in $\D(X)^\perp$ where $\phi\colon
S^1\to \R$ given by $\phi(V)=\|B_XV\|^2$ assumes its maximum
value. Differentiating  $\psi(t)= \phi(\cos t e_1+\sin t e_2)$
yields
$$
0=\psi'(0)=2\<B_Xe_1,B_Xe_2\>.
$$
Thus, there exist an orthonormal basis $\{\xi_1,\xi_2\}$ of
$T_z^\perp N$ and positive real numbers $\lambda_1, \lambda_2$
such that $B_Xe_r=\lambda_r\xi_r$ for $r=1,2$.
  We have from (\ref{eq:ge44}) that
\be\label{a8} \lambda_1\<\a(e_2,e_r),\xi_1\>=
\lambda_2\<\a(e_1,e_r),\xi_2\>,\;\;\; r=1,2. \ee Using that
(\ref{c}) holds for $E_1=E_2=X$, $V_1=e_s$ and $V_2=e_t$, $1\le
s,t\le 2$, we obtain
$$
\<\a(X,X),\a(e_r,e_r)\>=\<\a(X,e_r),\a(X,e_r)\> =\lambda_r^2,
$$
$$
\<\a(X,X),\a(e_1,e_2)\>=\<\a(X,e_1),\a(X,e_2)\>=0.
$$
Setting $\gamma_{st}^ r=\<\a(e_s,e_t),\xi_r\> =\gamma_{st}^
r,\;\a_r=\<\a(X,X),\xi_r\>\;\mbox{and}\;
D=\gamma_{11}^1\gamma_{22}^2-\gamma_{11}^2\gamma_{22}^1$, where
$1\le s,t\le 2$, it follows that \be\label{a11} a_1\gamma_{rr}^1 +
a_2\gamma_{rr}^2 = \lambda_r^2 \;\an\; a_1\gamma_{12}^1 +
a_2\gamma_{12}^2 = 0. \ee If we compute $Da_1$ and $Da_2$ from the
first equation in (\ref{a11}) and put the result into the second
we obtain an equation which  because of (\ref{a8}) is equivalent
to
$$
\lambda_1\lambda_2(\<\a(e_1,e_1),\a(e_2,e_2)\>
-\<\a(e_1,e_2),\a(e_1,e_2)\>)=0,
$$
and this gives (\ref{flat}) and concludes the proof of (\ref{c}).
\vspace{1ex}

   Because of Proposition \ref{prop:symm} it remains to show that
\be\label{eq:cyis} C(Y_1,Y_2,Y_3,Y_4)=0\;\;\; \mbox{for
all}\;\;Y_1,Y_2,Y_3,Y_4\in\H_z. \ee
  We divide the proof into the same two cases
considered before. \vspace{2ex}

\noindent {\it Case $1$.}  By Lemma \ref{prop:rankone} the linear
map $\H_z\to T_z^\perp N$, $\;Y\mapsto \a(Y,e)$ is surjective.
Hence $p\geq 2$, and we can take $X_1,X_2\in\H_z$ such that the
vectors $\xi_j:=\a(X_j,e)$, $j=1,2$, form an orthonormal normal
basis and $\a(Z,e)=0$ for all $Z\in \{X_1,X_2\}^\perp$. We have
from (\ref{eq:ge22}) and (\ref{eq:ge33}) that
$$
0=C(Z,X_j,e,E) = \<A_{\xi_j} Z,E\>
$$
for all $Z\in \{X_1,X_2\}^\perp$ and any $E\in T_zN$. Therefore,
$A_{\xi_j}Z =0$ for all $Z\in \{X_1,X_2\}^\perp$, that is,
$\{X_1,X_2\}^\perp\subset \Delta(z)$. Hence (\ref{eq:cyis}) holds
whenever $Y_i\in \{X_1,X_2\}^\perp$ for some $1\leq i\leq 4$.
Thus, it remains to prove that \be\label{ces}
C(X_1,X_2,X_2,X_1)=0. \ee Because of  (\ref{c}) we have
\be\label{eq:te} \<\a(X_i,X_j),\a(e,e)\>=\<\xi_i,\xi_j\>. \ee We
now may assume that $\a(e,e) = a\,\xi_2$ with $a\neq 0$.   We
obtain from (\ref{eq:te}) that \be\label{k3} \<A_{\xi_2}
X_1,X_2\>=0 \ee and
$$
a\<A_{\xi_2} X_j,X_j\> = 1,\;\;\; j=1,2.
$$
Since $a\neq 0$, it follows from the last equation that
\be\label{k4} \<A_{\xi_2} X_1,X_1\> = \<A_{\xi_2} X_2,X_2\>. \ee
In addition, (\ref{eq:ge22}) yields \be\label{k5}
0=C(X_i,X_j,e,X_i) = \<A_{\xi_j} X_i,X_i\> -\<A_{\xi_i}
X_i,X_j\>,\;\;\; i\neq j. \ee Then (\ref{ces}) follows from
(\ref{k3}), (\ref{k4}) and  (\ref{k5}), and the proof is completed
in this case. \vspace{2ex}

\noindent {\it Case $2$.}  Since (\ref{flat}) holds, it follows
from a result of E. Cartan (\cite{ca}; cf.\  Theorem~$1$ in
\cite{mo2}), that we may choose nonzero vectors $v_1,v_2$ in
$\D(X)^\perp$ (not necessarily orthogonal) such that
$\<\a(v_1,v_1),\a(v_2,v_2)\>=0$ and $\a(v_1,v_2)=0$. From
(\ref{c}) we have
$$
0=C(Z,v_1,v_2,Y)= - \<\a(Z,v_2),\a(Y,v_1)\>.
$$
Therefore, the subspaces $\a(\H_z,v_1), \a(\H_z,v_2)$ are
orthogonal lines spanned by $\eta_1=\a(X,v_1)$ and
$\eta_2=\a(X,v_2)$, respectively, which we may assume to have unit
length by rescaling $v_1$ and $v_2$ if necessary. In particular,
the kernel ${\mathcal H}_j$ of the linear map $F_j\colon\,\H_z\to
T_z^\perp N$ given by $F_j(Y)=\a(Y,v_j)$, $j=1,2$, has codimension
one in $\H_z$. On the other hand, for $Y\in {\mathcal H}_j$ and
any $Z\in \H_z$ we obtain using (\ref{eq:ge22}) that
$$
\<A_{\eta_j}Y,Z\>\!=\!\<\a(Y,Z),\a(X,v_j)\>
=\<\a(Y,v_j),\a(X,Z)\>=\<F_j(Y),\a( X,Z)\>=0.
$$
Let $Z_j$ be a unit vector in $\H_z$ orthogonal to ${\mathcal
H}_j$, $j=1,2$. Then
$$\begin{array}{l}C(Y_1,Y_2,Y_3,Y_4)=\<\a(Y_1,Y_4),
\a(Y_2,Y_3)\>-\<\a(Y_1,Y_3),\a(Y_2,Y_4)\>\vspace{1ex}\\\hspace{17ex}
=\sum_{j=1}^2(\<A_{\eta_j}Y_1,Y_4\>\<A_{\eta_j}Y_2,Y_3\>-
\<A_{\eta_j}Y_1,Y_3\>\<A_{\eta_j}Y_2,Y_4\>)\vspace{1ex}\\\hspace{17ex}
=\sum_{j=1}^2\<Y_1,Z_j\>\<Y_2,Z_j\>(\<A_{\eta_j}Z_j,Y_4\>\<A_{\eta_j}Z_j,Y_3\>-\vspace{1ex}
\\\hspace{20ex}\<A_{\eta_j}Z_j,Y_3\>\<A_{\eta_j}Z_j,Y_4\>)=0.\;\;\;\;\qed
\end{array}$$

Proposition \ref{casec} and the italicized text in its proof yield
the following for the case of Riemannian products.

\begin{corollary}\label{casecprod} Let
$f\colon\,L^p\times M^n\to\Q_c^{\,p+n+2}$ be an isometric
immersion of a Riemannian product. Assume that $f$ is of type $C$
at $z\in N^{p+n}=L^p\times M^n$. Then we have
\begin{itemize}
\item[$(i)$] If either $p\geq 3$ or $n\geq 3$ then $c=0$.
\item[$(ii)$] If $p+n\geq 3$ and $c=0$ then  $N^{p+n}$ is flat at $z$.
\end{itemize}
\end{corollary}

\end{document}